\documentclass[11pt,oneside,reqno]{amsart}
\usepackage{amssymb,amsmath,amsthm}
\usepackage{graphicx,color}
\usepackage{stackrel}
\usepackage{breqn}
\newcommand\numberthis{\addtocounter{equation}{1}\tag{\theequation}}

\usepackage{algorithm}
\usepackage[noend]{algpseudocode}

\usepackage{hyperref}
\theoremstyle{definition}

\numberwithin{equation}{section}
\usepackage[normalem]{ulem}
\newcommand{\txtr}{\textcolor{red}} 
\newcommand{\txtb}{\textcolor{blue}} 
\definecolor{DarkGreen}{rgb}{0,0.5,0.1} 

\newcommand\soutD{\bgroup\markoverwith
{\textcolor{DarkGreen}{\rule[.5ex]{2pt}{1pt}}}\ULon}
\newcommand\soutP{\bgroup\markoverwith
{\textcolor{blue}{\rule[.5ex]{2pt}{1pt}}}\ULon}
\newcommand{\Hm}[1]{\leavevmode{\marginpar{\tiny%
$\hbox to 0mm{\hspace*{-0.5mm}$\leftarrow$\hss}%
\vcenter{\vrule depth 0.1mm height 0.1mm width \the\marginparwidth}%
\hbox to
0mm{\hss$\rightarrow$\hspace*{-0.5mm}}$\\\relax\raggedright #1}}}

\begin{document}
%
\title[A well-conditioned MFS for Laplace equation]{
A well-conditioned Method of Fundamental Solutions for Laplace equation}
\author[P. R. S. Antunes]{
Pedro R.~S. Antunes}

\address{Sec\c{c}\~{a}o de Matem\'{a}tica, Departamento de Ci\^{e}ncias e Tecnologia, Universidade Aberta, Pal\'{a}cio Ceia, 1269-001 Lisboa and Group of Mathematical Physics, Faculdade de Ci\^{e}ncias da Universidade de Lisboa, Campo Grande, Edif\'{\i}cio C6
1749-016 Lisboa, Portugal} \email{prantunes@fc.ul.pt}

\date{\today}

\thanks{The research was partially supported by FCT, Portugal, through the scientific project UIDB/00208/2020.}

\begin{abstract}
The method of fundamental solutions (MFS) is a numerical method for solving boundary value problems involving linear partial differential equations. It is well-known that it can be very effective assuming regularity of the domain and boundary conditions. The main drawback of the MFS is that the matrices involved are typically ill-conditioned and this may prevent the method from achieving high accuracy.

In this work, we propose a new algorithm to remove the ill-conditioning of the classical MFS in the context of the Laplace equation defined in planar domains. The main idea is to expand the MFS basis functions in terms of harmonic polynomials. Then, using the singular value decomposition and Arnoldi orthogonalization we define well conditioned basis functions spanning the same functional space as the MFS's. Several numerical examples show that when possible to be applied, this approach is much superior to previous approaches, such as the classical MFS or the MFS-QR.
\end{abstract}

\maketitle
\section{Introduction}
%
The method of fundamental solutions (MFS) is a well known numerical method for the solution of boundary value problems for linear partial differential equations (PDEs) (\cite{Cheng,FK,Kupradze,mat-john}). It was introduced by Kupradze and Aleksidze \cite{Kupradze} and has been applied in the context of problems arising in acoustics (\cite{Alves-Antunes_2013,Antunes_drum,Barnett}), elasticity (\cite{Askour,Berger,Sarler}), fluid dynamics (\cite{Alves-Silvestre,kar-lesnic,Martins-Silvestre,Tsai}) or electromagnetism (\cite{Chen-electro,Young}).

It is a meshfree method for which the solution is approximated by a linear combination of shifts of the fundamental solution of the PDE to a set of points placed on an admissible source set. Since, by construction, the MFS approximation satisfies the PDE, we can focus just on the approximation of the boundary data, which is justified by density results (\cite{Alves,Bog}). The method may present a very fast convergence rate, under some regularity assumptions and its simplicity and effectiveness make it a very appealing numerical method. 

In the original formulation the application of the MFS was restricted to the numerical solution of homogeneous linear PDEs defined in smooth domains. Recent studies have extended the range of applications of the MFS, which includes the application to non-homogeneous PDEs (\cite{Alves-Valtchev,chennh}), to non smooth boundaries (\cite{Alves-Leitao,Antunes_Valtchev}) and non smooth boundary conditions (\cite{Alves-Martins-Valtchev}). There are two major issues in the MFS yet to be resolved: the choice of the source points and the ill-conditioning. The optimal location of the source points has been widely addressed in the literature (\cite{Alves,Alves-Antunes_2005,Alves-Antunes_2013,Barnett,Chen-Kara,Cheng,zcli,Hon,kar1,Li-CS-Kara}) and several choices have been advocated to be effective. Some works have proposed techniques to alleviate the ill-conditioning of the MFS \cite{MFSQR,Ant2,CSChen,raman}, but none of these approaches seem to completely solve the problem.  

In this work we propose a new technique to remove the ill-conditioning in the MFS. The main idea is to consider an expansion of the MFS basis functions in terms of harmonic polynomials. Then, using the singular value decomposition and Arnoldi orthogonalization we define well conditioned basis functions spanning the same functional space as the original MFS basis functions. To the best of our knowledge, this is the first approach able to remove the ill-conditioning of the classical MFS in the context of general planar domains.

\section{The classical MFS}
\label{classicalmfs}
Let $\Omega$ be a bounded and smooth planar domain. We will consider the numerical solution of the Dirichlet boundary value problem involving the Laplace equation,
\begin{equation}
\label{bvp} \left\{ \begin{array}{cl}
\Delta u=0 & \mbox{ in }\Omega,\\
u=g & \mbox{ on }\partial\Omega,\end{array}\right.\end{equation}
for a given function $g$ defined on $\partial\Omega$.

As in \cite{MFSQR}, we will use the terminology Direct-MFS to refer to the classical MFS approach that we briefly describe in this section.

Let $\Phi$ be a fundamental solution of the operator $-\Delta$,
\[\Phi(x)=-\frac{1}{2\pi}\log\left|x\right|.\]
This fundamental solution is analytic, except at the origin, where it has a singularity. The Direct-MFS approximation is a linear combination
\begin{equation}
\label{direct-mfs}
u^{\text{Direct}}_N(x)=\sum_{n=1}^Nc^{\text{Direct}}_n\Phi(x-y_n),
\end{equation}
where each base function is a translation of the fundamental solution to some source point $y_n$ placed on some admissible source set $\hat{\Gamma}$ that does not intersect $\bar{\Omega}$. 

The approximation can be mathematically justified by density results stating that
\[\text{span}\left\{\Phi(\bullet-y)|_{\Omega}:y\in\hat{\Gamma}\right\}\]
is dense in $L^2(\partial\Omega)$ \cite{Alves,Bog}. This Direct-MFS linear combination can be augmented with some extra basis functions. For example for ensuring completeness it is necessary to add the constant function \cite{Alves,JTChen2} and to improve the MFS accuracy in the context of non smooth problems, some singular particular solutions may be added \cite{Alves-Leitao,Alves-Martins-Valtchev}.

By construction, the MFS approximation \eqref{direct-mfs} satisfies the Laplace equation in $\Omega$ and the coefficients can be calculated by forcing the boundary conditions at some boundary points. We consider $M$ (with $M\geq N$) collocation points $x_i$, $i=1,2,...,M$, and solve
\begin{equation}
\label{sistema-direct}
\mathbf{A}^{\text{Direct}} \cdot \mathbf{c}^{\text{Direct}}=\mathbf{G},
\end{equation}
where
\[\mathbf{A}^{\text{Direct}}=\left[\Phi(x_i-y_j)\right]_{M\times N},\quad\mathbf{G}=\left[g(x_i)\right]_{M\times 1}\]
and $\mathbf{c}^{\text{Direct}}$ is a vector containing all the coefficients of the Direct-MFS linear combination \eqref{direct-mfs}. In all computations in this work we shall consider oversampling and \eqref{sistema-direct}
will be solved in the least-squares sense. We will fix $M=2N$ but the method is not much sensitive to other choices of $M>N$.

The Direct-MFS may provide highly accurate approximations, even considering just a few source points. For instance, assume that $\Omega$ is the unitary disk and $g$ is the trace of an entire function. If we define $M=N$ collocation points uniformly distributed on $\partial\Omega$ and the source points are placed on a circle of radius equal to $R$, we have exponential convergence (cf.~\cite{Katsurada1,Katsurada2}) and there exists a constant $C>0$ that does not depend on $N$, such that
\begin{equation}
\label{expconv}
\left\|g-u_N\right\|_{L^2(\partial\Omega)}\leq CR^{-N}.
\end{equation}
Moreover, this bound might suggest that the convergence is faster for large values of $R$, that is, choosing the artificial boundary, where we place the source points, far from the boundary. However, under similar assumptions we know that the condition number also increases exponentially (\cite{Kitagawa,smy-kar}), 
\begin{equation}
\label{condnumdirect}
cond_2\left(\mathbf{A}^{\text{Direct}}\right)\sim\frac{1}{2}\log(R)R^{\frac{N}{2}}.
\end{equation}
Therefore, for large values of $R$, the Direct-MFS is highly ill-conditioned, which affects the accuracy and prevents the exponential convergence to be achieved.

In the following section we will describe a different approach that reduces the problem of ill-conditioning of the Direct-MFS.

\section{The MFS-QR}
The MFS-QR was introduced in~\cite{MFSQR} as a technique to reduce the ill-conditioning of the Direct-MFS. The source points are assumed to be placed on a circle of radius $\frac{1}{\epsilon}$, 
\begin{equation}
\label{sourcepoints}
y_j=\frac{1}{\epsilon}\left(\cos(\alpha_j),\sin(\alpha_j)\right),\ j=1,...,N,\quad\alpha_j=\frac{2\pi j}{N},
\end{equation}
where $\epsilon$ is chosen such that $\frac{1}{\epsilon}>R_\Omega:=\max_{x\in\partial\Omega}\left\|x\right\|.$

Each source point can be written as $\frac{1}{\epsilon}\left(\cos(\alpha),\sin(\alpha)\right)$,
for some $\alpha\in[0,2\pi[$ and dropping the constant $-1/(2\pi)$, we have the following expansion of the corresponding MFS base function in polar coordinates (\cite{MFSQR,JTChen}),
\begin{align*}
\psi(r,\theta)&=\log\left(\sqrt{\left(r\cos(\theta)-\frac{1}{\epsilon}\cos(\alpha)\right)^2+\left(r\sin(\theta)-\frac{1}{\epsilon}\sin(\alpha)\right)^2}\right)\\
&=\log\left(\sqrt{r^2+\frac{1}{\epsilon^2}-\frac{2r}{\epsilon}\cos(\theta-\alpha)}\right)\\
&=\log\left(\frac{1}{\epsilon}\sqrt{\epsilon^2r^2+1-2r\epsilon\cos(\theta-\alpha)}\right)\\
&=-\log(\epsilon)+\log\left(\sqrt{\epsilon^2r^2+1-2r\epsilon\cos(\theta-\alpha)}\right)\\
&=-\log(\epsilon)+\frac{1}{2}\log\left(\epsilon^2r^2+1-2r\epsilon\cos(\theta-\alpha)\right)\\
&=-\log(\epsilon)+\frac{1}{2}\log\left(\epsilon^2r^2+1-r\epsilon\left(e^{i(\theta-\alpha)}+e^{-i(\theta-\alpha)}\right)\right)\\
&=-\log(\epsilon)+\frac{1}{2}\log\left[\left(1-\epsilon re^{i(\theta-\alpha)}\right)\left(1-\epsilon re^{-i(\theta-\alpha)}\right)\right]\\
&=-\log(\epsilon)+\frac{1}{2}\left[\log\left(1-\epsilon re^{i(\theta-\alpha)}\right)+\log\left(1-\epsilon re^{-i(\theta-\alpha)}\right)\right]\\
&=-\log(\epsilon)+\frac{1}{2}\left[\sum_{m=1}^\infty\frac{-1}{m}\left(\epsilon r e^{i(\theta-\alpha)}\right)^m+\sum_{m=1}^\infty\frac{-1}{m}\left(\epsilon r e^{-i(\theta-\alpha)}\right)^m\right]\\
&=-\log(\epsilon)-\sum_{m=1}^{\infty}\frac{r^m\epsilon^m}{m}\frac{e^{im(\theta-\alpha)}+e^{-im(\theta-\alpha)}}{2}\numberthis\label{expcomp} \\
&=-\log(\epsilon)-\sum_{m=1}^{\infty}\frac{r^m\epsilon^m}{m}\cos(m(\theta-\alpha))\\
&=-\log(\epsilon)-\sum_{m=1}^{\infty}\frac{r^m\epsilon^m}{m}\left(\cos(m\theta)\cos(m\alpha)+\sin(m\theta)\sin(m\alpha)\right).\numberthis \label{expansao}
\end{align*}
Therefore, the MFS basis functions $\psi_j(r,\theta)$ associated to $N$ source points are given by
\begin{small}
\[\begin{bmatrix}
    \psi_1(r,\theta) \\
     \psi_2(r,\theta) \\
    \vdots  \\
   \psi_N(r,\theta)\end{bmatrix}=\begin{bmatrix}
    -1 & -\cos(\alpha_1) & -\sin(\alpha_1) & -\cos(2\alpha_1) & -\sin(2\alpha_1)& -\cos(3\alpha_1) & -\sin(3\alpha_1) & \dots \\
    -1 & -\cos(\alpha_2)  & -\sin(\alpha_2) & -\cos(2\alpha_2)  & -\sin(2\alpha_2)& -\cos(3\alpha_2)  & -\sin(3\alpha_2)  & \dots   \\
    \vdots & \vdots & \vdots& \vdots & \vdots& \vdots & \vdots & \ddots  \\
    -1 & -\cos(\alpha_N) & -\sin(\alpha_N) & -\cos(2\alpha_N) & -\sin(2\alpha_N)& -\cos(3\alpha_N) & -\sin(3\alpha_N)  & \dots
\end{bmatrix}\]
\begin{equation}
\label{expansao1}
\begin{bmatrix}
    \log(\epsilon) & & &  & & &  & \   \\
    & \epsilon &  & & & &  &   \\
     & & \epsilon &  && & &  \\
     & & & \frac{\epsilon^2}{2}   && & &  \\
          & & &   &\frac{\epsilon^2}{2} & & &  \\
               & & & & & \frac{\epsilon^3}{3}    & &  \\
                              & & & & & & \frac{\epsilon^3}{3}     &  \\
   &  & &  & &  & & \ddots 
\end{bmatrix}\begin{bmatrix}
    1  \\
    r\cos(\theta)  \\
    r\sin(\theta)  \\
        r^2\cos(2\theta)  \\
    r^2\sin(2\theta)  \\
        r^3\cos(3\theta)  \\
    r^3\sin(3\theta)  \\
    \vdots \\
\end{bmatrix}
\end{equation}
\end{small}
and after truncating this expansion, considering the sum in \eqref{expansao} just up to $m=p$ (such that $2p+1>N$) we obtain a factorization
\begin{equation}
\label{vect}
\mathbf{\Theta}(r,\theta) =\mathbf{B}\cdot\mathbf{D}\cdot\mathbf{F}(r,\theta).
\end{equation}
The MFS-QR involves the calculation of a $\textbf{QR}$ factorization of the matrix $\mathbf{B}$,
\[\mathbf{B}= \mathbf{Q}\cdot\mathbf{R},\]
and we define a new set of functions spanning the same functional space as Direct-MFS basis functions given by
\[\mathbf{\Psi}(r,\theta)=\tilde{\mathbf{R}}\cdot\mathbf{F}(r,\theta),\]
defining
\[\tilde{\mathbf{R}}=\tilde{\mathbf{T}}\circ \mathbf{R},\]
where
\[\hspace{-1cm}\tilde{\mathbf{T}}=\begin{bmatrix}
  1& \epsilon/ \log(\epsilon) & \epsilon/ \log(\epsilon) & \epsilon^2/(2 \log(\epsilon))  & \epsilon^2/ (2\log(\epsilon))& \dots &  \epsilon^p/(p \log(\epsilon)) & \epsilon^p/( \log(\epsilon)p)\\
    & 1 & 1&\epsilon/2 & \epsilon/2 & \dots & \epsilon^{p-1}/p& \epsilon^{p-1}/p \\
    &  & 1&\epsilon/2 & \epsilon/2 & \dots & \epsilon^{p-1}/p& \epsilon^{p-1}/p \\
    &  & &1 & 1 & \dots & \epsilon^{p-2}/p& \epsilon^{p-2}/p \\
    &  & & & 1 & \dots & \epsilon^{p-2}/p& \epsilon^{p-2}/p \\
   &  & &  & &  \ddots  & & \ddots 
\end{bmatrix}_{N\times (2p+1)}\]
and $\circ$ denotes the Hadamard product of matrices,
\[(\mathbf{A}\circ \mathbf{B})_{i,j}=(\mathbf{A})_{i,j}.(\mathbf{B})_{i,j}.\]
The MFS-QR approximation is a linear combination 
\begin{equation}
\label{mfs-qr}
u^{QR}_N(r,\theta)=\sum_{n=1}^Nc^{QR}_n\mathbf{\Psi}_n(r,\theta)
\end{equation}
and the solution of the boundary value problem \eqref{bvp} is obtained through the solution of the linear system
\begin{equation}
\label{sistema-qr}
\mathbf{A}^{QR}.\mathbf{c}^{QR}=\mathbf{G},
\end{equation}
where $\mathbf{A}^{QR}=\left[\mathbf{\Psi}(x_i)\right]^T$ or equivalently
\[\left[\mathbf{F}(x_i)\right]^T.\tilde{\mathbf{R}}^T.\mathbf{c}^{QR}=\mathbf{G}.\]

As was reported in \cite{MFSQR}, the MFS-QR removes completely the ill-conditioning when $\partial\Omega$ and the artificial boundary are concentric circles. However, if the domain is not a disk, though at a lower rate than the Direct-MFS, the condition number also grows exponentially as the number of basis functions increase and this can be prohibitive for large values of $N$. The main reason for this growth is the increasing powers of $r$ in the harmonic polynomials in $\mathbf{F}(r,\theta)$ that are highly ill-conditioned.

Moreover, by construction, the MFS-QR assumes the artificial boundary to be a circle and this may be too restrictive when applying the MFS to some geometries, such as domains that are elongated in one direction or domains with re-entrant corners for which, in general, placing the source points on a circle may not be a suitable choice.

Next, we will introduce a new technique to remove the ill-conditioning of the MFS in the context of general planar domains and artificial boundaries, provided the latter satisfy a geometric constraint defined in next section.

\section{A new technique to remove the ill-conditioning of the Direct-MFS}
\label{neqtech}
\begin{figure}[ht]
\centering \includegraphics[width=0.38\textwidth]{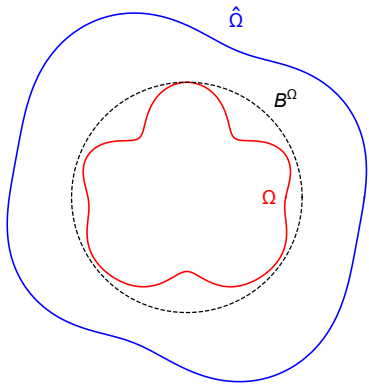}
\caption{The setting for the numerical approach. We assume that there exists a ball $B^\Omega$ such that $\Omega\subset B^\Omega\subset \hat{\Omega}.$} \label{fig:setting}
\end{figure}

In this section, we will describe a new technique to remove the ill-conditioning of the Direct-MFS. We assume that $\Omega$ is a smooth star shaped planar domain and for simplicity, we assume that it contains the origin. We will denote by $\hat{\Omega}$ an open set for which $\bar{\Omega}\subset\hat{\Omega}$ and will take the boundary of $\hat{\Omega}$ as an artificial curve for the MFS. We assume that $\hat{\Omega}$ is chosen in such a way that there exists a ball $B^\Omega$ for which we have $\Omega\subset B^\Omega\subset \hat{\Omega},$ as illustrated in Figure~\ref{fig:setting}.
We will denote by $(r,\theta)$ the polar coordinates of an arbitrary point $x\in\mathbb{R}^2$, that is
\[x=r(\cos(\theta),\sin(\theta))\]
and given a sample of $N$ source points placed outside $\Omega$ we calculate their polar coordinates, which we assume to be of the form $\left(\frac{1}{\epsilon_j},\alpha_j\right)$ where $\alpha_j\in[0,2\pi[$, for $j=1,...,N$. Thus, using the fact that for any $z\in\mathbb{C}$, we have
\[\cos(z)=\frac{e^{iz}+e^{-iz}}{2},\] 
the MFS basis functions corresponding to a source point $y_j$ have the following expansion, directly obtained from \eqref{expcomp}, 
\begin{align*}\psi_j(r,\theta)&=-\log(\epsilon_j)-\sum_{m=1}^{\infty}\frac{r^m\epsilon_j^m}{m}\frac{e^{im(\theta-\alpha)}+e^{-im(\theta-\alpha)}}{2}\\
&=-\log(\epsilon_j)-\sum_{m=1}^{\infty}\frac{r^m\epsilon_j^m}{m}\frac{e^{im\theta}e^{-im\alpha}+e^{-im\theta}e^{im\alpha}}{2}\\
&=-\log(\epsilon_j)-\sum_{m=1}^{\infty}\left(\frac{r}{R_\Omega}\right)^m\left(\epsilon_j R_\Omega\right)^m\frac{e^{im\theta}e^{-im\alpha}+e^{-im\theta}e^{im\alpha}}{2m}\numberthis \label{eqnultima}
\end{align*}
where for convenience we multiplied and divided each term in the sum by $R_\Omega^m.$ Note that, by the definition of $R_\Omega,$ for all points in $\bar{\Omega}$ we have $\frac{r}{R_\Omega}\leq1$ and the assumption that there exists a ball $B^\Omega$ for which we have $\Omega\subset B^\Omega\subset \hat{\Omega}$
ensures that
\[\max_{j=1,...,N}\left(\epsilon_jR_\Omega\right)<1\]
which implies convergence of the series.
Formally, we can expand all the MFS basis functions defined for $(r,\theta)\in\bar{\Omega}$ in terms of harmonic polynomials as follows
\[\boldsymbol{\psi}(r,\theta)=\mathbf{M_\infty}\cdot \boldsymbol{F_\infty}(r,\theta)\]
   where
   \[\boldsymbol{\psi}(r,\theta)=\begin{bmatrix}
    \psi_1(r,\theta) &
     \psi_2(r,\theta) &
    \dots  &
   \psi_N(r,\theta)\end{bmatrix}^T,\]
   \[\mathbf{M_\infty}=\begin{bmatrix}
    -\log(\epsilon_1) & -\frac{(\epsilon_1R_\Omega)e^{-i\alpha_1}}{2} & -\frac{(\epsilon_1R_\Omega)^2e^{-2i\alpha_1}}{4} & \dots & -\frac{(\epsilon_1R_\Omega)e^{i\alpha_1}}{2}& -\frac{(\epsilon_1R_\Omega)^2e^{2i\alpha_1}}{4} & \dots  \\
    -\log(\epsilon_2) & -\frac{(\epsilon_2R_\Omega)e^{-i\alpha_2}}{2}  & -\frac{(\epsilon_2R_\Omega)^2e^{-2i\alpha_2}}{4} & \dots  & -\frac{(\epsilon_2R_\Omega)e^{i\alpha_2}}{2}& -\frac{(\epsilon_2R_\Omega)^2e^{2i\alpha_2}}{4}  & \dots   \\
    \vdots & \vdots & \vdots & \ddots& \vdots & \vdots & \ddots  \\
    -\log(\epsilon_N) & -\frac{(\epsilon_NR_\Omega)e^{-i\alpha_N}}{2} & -\frac{(\epsilon_NR_\Omega)^2e^{-2i\alpha_N}}{4} & \dots & \frac{(\epsilon_NR_\Omega)e^{i\alpha_N}}{2}& \frac{(\epsilon_NR_\Omega)^2e^{2i\alpha_N}}{4} &  \dots
\end{bmatrix}\]
and
\[\boldsymbol{F_\infty}(r,\theta)=\begin{bmatrix}
1 &
\frac{r}{R_\Omega}e^{i\theta}&
\left(\frac{r}{R_\Omega}\right)^2e^{2i\theta} &
\dots &  
\frac{r}{R_\Omega}e^{-i\theta} &
\left(\frac{r}{R_\Omega}\right)^2e^{-2i\theta} &
\dots
\end{bmatrix}^T.\]
Note that since $R_\Omega$ is just a constant,  the (infinite) matrix $\mathbf{M}$ depends just on the source points while $\boldsymbol{F}(r,\theta)$ depends just on $\bar{\Omega}.$ The expansions can be truncated in such a way that the residual is smaller than machine precision, for instance considering the expansion \eqref{eqnultima} just up to $m=p_0,$ provided $2p_0+1\geq N$. Actually, we can determine $p_0$ taking into account that
\begin{align*}\psi_j(r,\theta)&=-\log(\epsilon_j)-\sum_{m=1}^{\infty}\left(\frac{r}{R_\Omega}\right)^m\left(\epsilon_j R_\Omega\right)^m\frac{e^{im\theta}e^{-im\alpha}+e^{-im\theta}e^{im\alpha}}{2m}\\
&=-\log(\epsilon_j)-\sum_{m=1}^{p_0}\left(\frac{r}{R_\Omega}\right)^m\left(\epsilon_j R_\Omega\right)^m\frac{e^{im\theta}e^{-im\alpha}+e^{-im\theta}e^{im\alpha}}{2m}-\\
&\sum_{m=p_0+1}^{\infty}\left(\frac{r}{R_\Omega}\right)^m\left(\epsilon_j R_\Omega\right)^m\frac{e^{im\theta}e^{-im\alpha}+e^{-im\theta}e^{im\alpha}}{2m}
\end{align*}
and we can estimate the residual by
\begin{align*}
\left|\sum_{m=p_0+1}^{\infty}\left(\frac{r}{R_\Omega}\right)^m\left(\epsilon_j R_\Omega\right)^m\frac{e^{im\theta}e^{-im\alpha}+e^{-im\theta}e^{im\alpha}}{2m}\right|\leq\\
\sum_{m=p_0+1}^{\infty}\left|\left(\frac{r}{R_\Omega}\right)^m\left(\epsilon_j R_\Omega\right)^m\frac{e^{im\theta}e^{-im\alpha}+e^{-im\theta}e^{im\alpha}}{2m}\right|\leq\\
\sum_{m=p_0+1}^{\infty}\left|\frac{\left(\max_j\epsilon_j\ R_\Omega\right)^m}{m}\right|=\left(\max_j\epsilon_j\ R_\Omega\right)^{p_0+1}\Phi^{HL}(\max_j\epsilon_j\ R_\Omega,1,p_0+1),
\end{align*}
where $\Phi^{HL}$ is the Hurwitz-Lerch transcendent function. Thus, we define $p_0$ to be the smallest integer such that
\[\left(\max_j\epsilon_j\ R_\Omega\right)^{p_0+1}\Phi^{HL}(\max_j\epsilon_j\ R_\Omega,1,p_0+1)\leq\epsilon.\]
After determining $p:=\max(p_0,\lceil \frac{N-1}{2} \rceil)$ we truncate the expansion and obtain
\begin{equation}
    \label{expdisc}
\boldsymbol{\psi}(r,\theta)=\mathbf{M}\cdot \boldsymbol{F}(r,\theta),
   \end{equation}
where $\mathbf{M}$ is a $N\times(2p+1)$ matrix and $\boldsymbol{F}(r,\theta),$ a vector-valued function with $2p+1$ components that are the truncated versions of $\mathbf{M_\infty}$ and $\boldsymbol{F_\infty}(r,\theta),$ respectively. Note that there is a slight abuse of notation in \eqref{expdisc} in the sense that the equality shall be understood as an approximation with accuracy at machine precision level.

We will propose a new technique to remove the ill-conditioning in the Direct-MFS by performing a suitable change of basis. For this purpose, an important remark is that in \eqref{expdisc} we can multiply from the left by any non singular matrix. This corresponds to applying a change of the basis functions without modifying the functional space spanned by them. Thus, an interesting question is to identify a suitable matrix to be applied from the left in \eqref{expdisc}.

Given a set of collocation points and their polar coordinates 
\[x_j=(r_j,\theta_j),\ j=1,...,P,\] 
the (transpose of the) matrix of the Direct-MFS system \eqref{sistema-direct}, $\mathbf{A}^{Direct}$, is obtained directly from evaluating \eqref{expdisc} at the collocation points. At this point we can identify two possible sources for the ill-conditioning of the Direct-MFS:
\begin{itemize}
\item The matrix $\mathbf{M}$ is ill-conditioned 
\item the increasing powers of $\frac{r}{R_\Omega}$ in $\mathbf{F}(r,\theta)$ also generate ill-conditioning.
\end{itemize}

We will propose a technique to tackle these two sources of ill-conditioning respectively in next two sections.

\subsection{Reducing the effect of the ill-conditioning of matrix $\mathbf{M}$}
\label{reduce}
For now, we will assume that we can write the functions in $\mathbf{F}(r,\theta)$ in terms of a well conditioned basis,
\begin{equation}
\label{expbc}
\mathbf{F}(r,\theta)=\mathbf{K}\cdot\mathbf{J}(r,\theta),
\end{equation}
where $K$ is an invertible matrix and $\mathbf{J}(r,\theta)$ is a vector-valued function built by a set of well conditioned basis functions that are particular solutions of the Laplace equation. Note that from \eqref{expdisc} we can write

\begin{equation}
    \label{expdisc2}
\boldsymbol{\psi}(r,\theta)=\underbrace{\mathbf{M}\cdot \mathbf{K}}_{:=\mathbf{M}_1}\cdot\mathbf{J}(r,\theta):=\mathbf{M}_1\cdot\mathbf{J}(r,\theta).
   \end{equation}

We start by calculating the singular value decomposition of matrix $\mathbf{M}_1.$ This allows us to obtain the factorization
\[\mathbf{M}_1=\mathbf{U}\cdot\mathbf{S}\cdot\mathbf{V}^\ast,\]
where $\mathbf{U}$ and $\mathbf{V}$ are unitary and $\mathbf{S}$ is diagonal with non negative entries. Thus, multiplying by the matrix $\mathbf{U^\ast}$ from the left we obtain
\begin{equation}
    \label{contasmatrizes}
\mathbf{U^\ast}\cdot\mathbf{M}_1=\underbrace{\mathbf{U^\ast}\cdot\mathbf{U}}_{I}\cdot\mathbf{S}\cdot\mathbf{V}^\ast=\mathbf{S}\cdot\mathbf{V}^\ast.\end{equation}
We know that $\mathbf{S}$ has the same dimensions as $\mathbf{M}_1$ and since we took $2p+1\geq N$ we can write
\[\mathbf{S}=\left[\begin{array}{c|c}
\mathbf{S}_1&\mathbf{0} 
\end{array}
\right],\]
where $\mathbf{S}_1$ is a diagonal square matrix and $\mathbf{0}$ denotes a block matrix with all entries equal to zero and we have
\[\mathbf{S}\cdot\mathbf{V}^\ast=\left[\begin{array}{c|c}
\mathbf{S}_1&\mathbf{0} 
\end{array}
\right]\left[\begin{array}{c}
\mathbf{V}^\ast_1\\
\hline \mathbf{V}^\ast_2
\end{array}
\right]=\mathbf{S}_1\cdot\mathbf{V}^\ast_1+\mathbf{0}\cdot\mathbf{V}^\ast_2=\mathbf{S}_1\cdot\mathbf{V}^\ast_1,\]
where $\mathbf{V}^\ast_1$ is the matrix composed by the first $N$ rows of $\mathbf{V}^\ast$. Now, multiplying from the left by the matrix $\mathbf{S}_1^{-1}$ in \eqref{contasmatrizes} we obtain
\[\mathbf{S}_1^{-1}\cdot\mathbf{U^\ast}\cdot\mathbf{M}_1=\mathbf{S}_1^{-1}\cdot\mathbf{S}\cdot\mathbf{V}^\ast=\underbrace{\mathbf{S}_1^{-1}\cdot\mathbf{S}_1}_{I}\cdot\mathbf{V}^\ast_1=\mathbf{V}^\ast_1.\]

Thus, we obtain a new set of basis functions directly from \eqref{expdisc} (formally) multiplying from the left by the matrix $\left(\mathbf{S}_1^{-1}\cdot\mathbf{U^\ast}\right),$
\[\boldsymbol{\phi}(r,\theta)=\begin{bmatrix}
    \phi_1(r,\theta) &\phi_2(r,\theta) &\dots  &\phi_N(r,\theta)\end{bmatrix}^T\]
    defined by
\begin{equation}
    \label{expdisc3}
\boldsymbol{\phi}(r,\theta)=\left(\mathbf{S}_1^{-1}\cdot\mathbf{U^\ast}\right)\boldsymbol{\psi}(r,\theta)=\mathbf{V}^\ast_1\cdot \mathbf{J}(r,\theta).
   \end{equation}
Note that from \eqref{expdisc3} instead of the multiplication by the matrix $\left(\mathbf{S}_1^{-1}\cdot\mathbf{U^\ast}\right)$ to perform the change of basis, for practical purposes we simply calculate the product $\mathbf{V}^\ast_1\cdot \mathbf{J}(r,\theta)$ in order to have the new set of basis functions. 

We will call MFS-SVD to the approximation made by a linear combination of the new basis functions of $\boldsymbol{\phi}(r,\theta),$
\begin{equation}
\label{mfs-svd}
u_N^{\text{SVD}}(r,\theta)=\sum_{n=1}^Nc^{\text{SVD}}_n\phi_n(r,\theta)
\end{equation}
and the coefficients can be determined by solving
\begin{equation}
\label{sistema-novo}
[\boldsymbol{\phi}(x_i)]^T .\ \mathbf{C}^{\text{SVD}}=\mathbf{G}.
\end{equation}

Each component of the vector valued function $\mathbf{J}(r,\theta)$ satisfies the Laplace equation. Thus, each component of $\boldsymbol{\phi}(r,\theta)$, which is a linear combination of the components of $\mathbf{J}(r,\theta)$ is also a particular solution of the Laplace equation. Thus, since $u$ and $u^{SVD}_N$ are both harmonic functions, by the maximum principle we have 
\begin{equation}
\label{maxprinc}
\left\|u-u^{SVD}_N\right\|_{L^\infty(\Omega)}\leq\left\|u-u^{SVD}_N\right\|_{L^\infty(\partial\Omega)}=\left\|g-u^{SVD}_N\right\|_{L^\infty(\partial\Omega)}.
\end{equation}
In all the numerical simulations we will measure the error by evaluating the error $\left\|g-u^{SVD}_N\right\|_{L^\infty(\partial\Omega)}$ at 10001 boundary points.

\begin{figure}[ht]
\centering 
\includegraphics[width=0.48\textwidth]{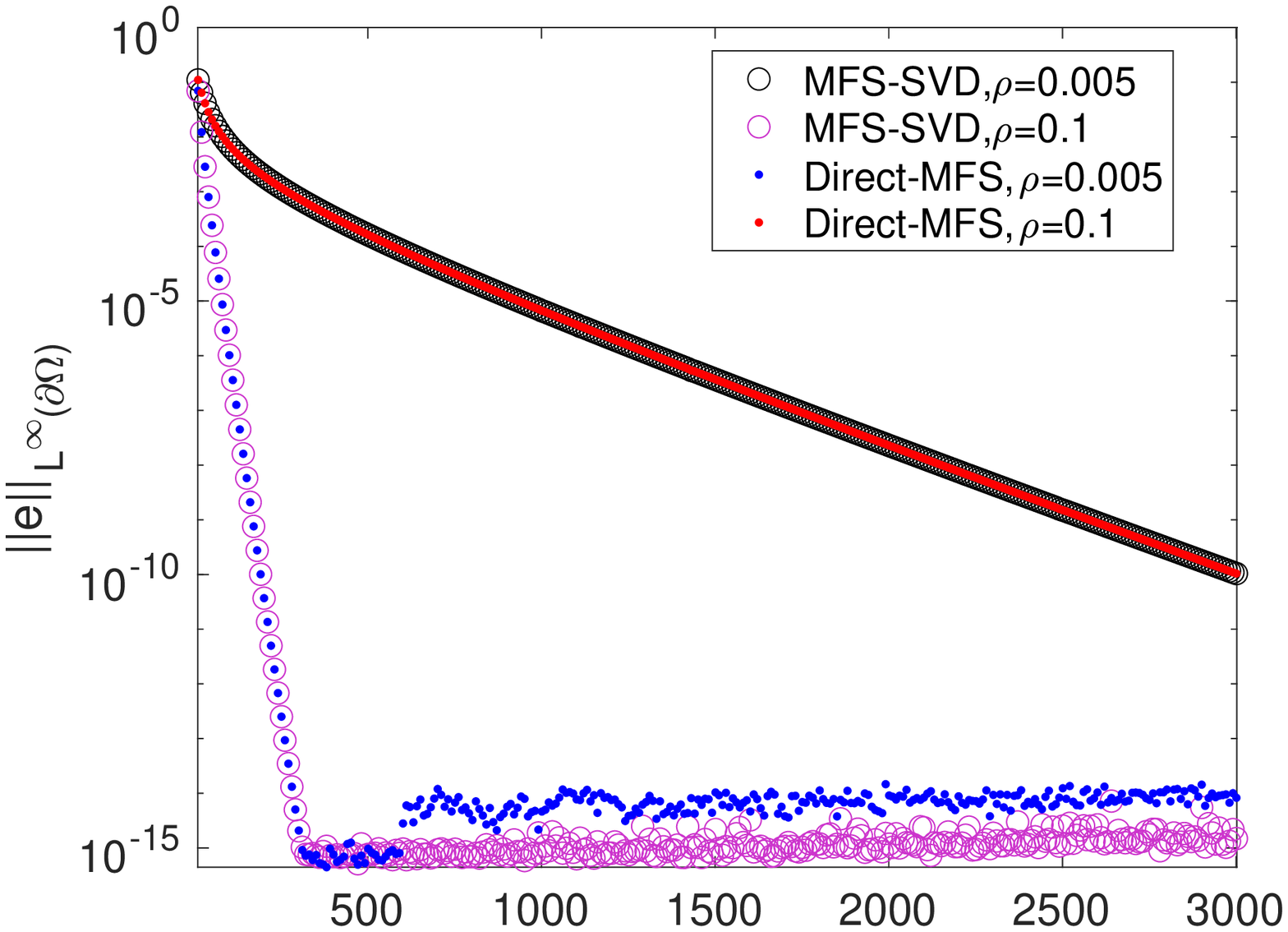}
\includegraphics[width=0.48\textwidth]{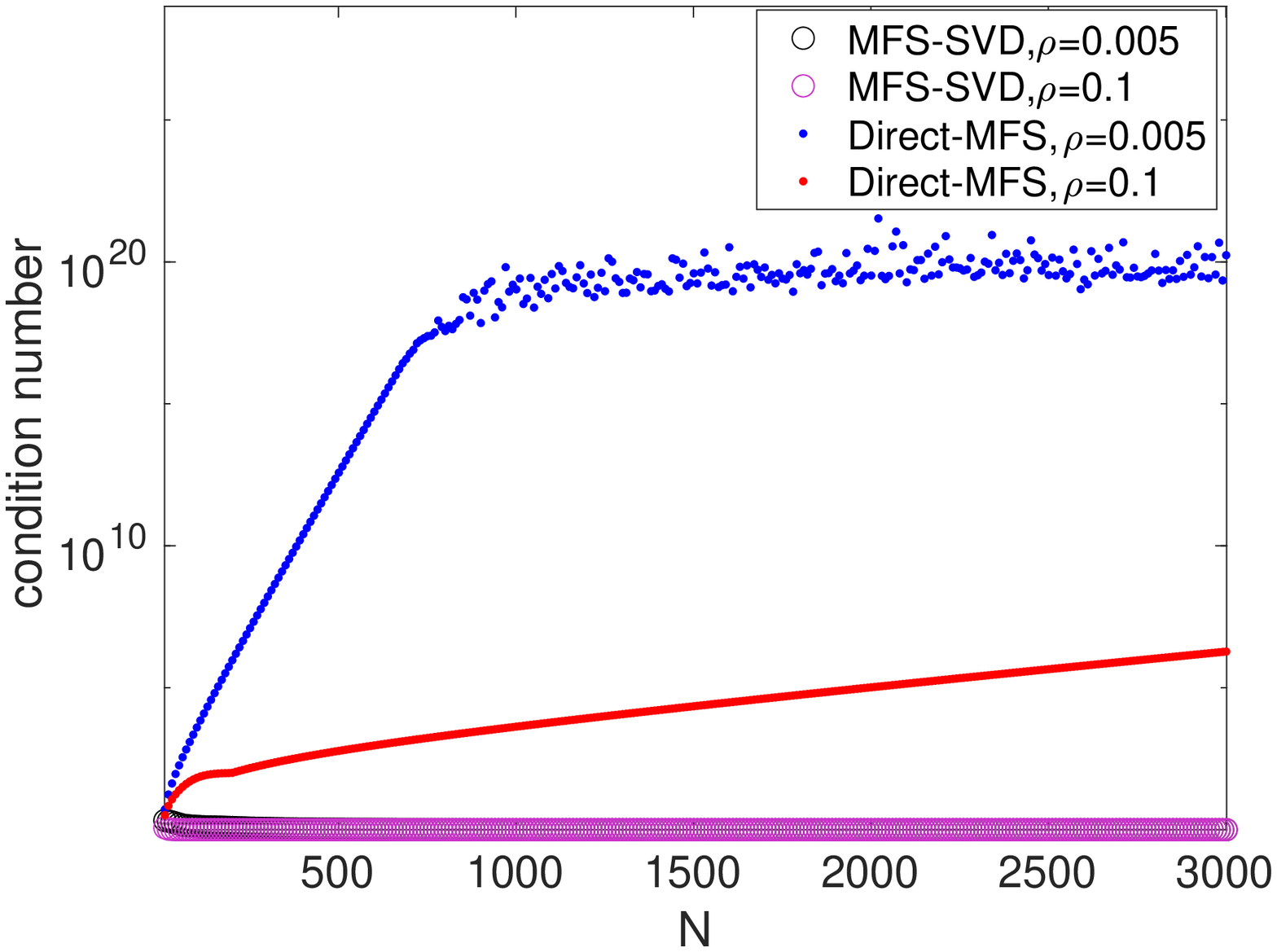}
\caption{Plot of the $L^\infty$ norm of the error of the approximations obtained with Direct-MFS and MFS-SVD, as a function of $N$ (left plot) and plot of the condition number of the system matrices of the two methods (right plot).} \label{fig:disco_resultados}
\end{figure}

To illustrate the numerical technique to remove the ill-conditioning of matrix $\mathbf{M}$ developed in this section we consider the unit disk. In the first numerical simulation we choose the artificial boundary to be the circle of radius equal to $1+\rho$, for a positive parameter $\rho$ and take the boundary data given by $g(x,y)=x^2y^3.$ Previous studies (eg.~\cite{Chen-Kara}) indicate that the source points should be chosen close to the 
boundary for achieving an optimal rate of convergence of the MFS for boundary value problems with non-harmonic boundary conditions. Figure~\ref{fig:disco_resultados}-left shows the plot of $L^\infty$ norm of the error of the approximations obtained with Direct-MFS and MFS-SVD, for $\rho=0.005$ and $\rho=0.1$. The right plot of the same figure shows the condition number of the system matrices. We can observe that when $\rho=0.005$ with obtain similar errors for Direct-MFS and MFS-SVD. The condition number of Direct-MFS grows exponentially, while the condition number of the MFS-SVD is always of order O(1), independently of the value of $N$. Taking $\rho=0.1$, the convergence is faster. The Direct-MFS and MFS-SVD have similar accuracy for $N<N_0\approx 600$. For $N>N_0$, the condition number of the Direct-MFS blows up and the two approaches present slightly different results.

Note that the components of $\mathbf{F}(r,\theta)$ are orthogonal functions in the unit disk, which means that in this case we can assume that the matrix $\mathbf{K}$ in \eqref{expbc} is simply the identity matrix. Thus, we can apply the technique developed in section~\ref{reduce} for arbitrary artificial boundary. To illustrate this fact, we consider the case where $\hat{\Omega}$ is the domain with boundary defined by
\[\Gamma=\left\{\left(4\gamma(t)\cos(t)-1,4\gamma(t)\sin(t)-1\right)\in\mathbb{R}^2:\ t\in[0,2\pi[\right\},\]
where $\gamma(t)=e^{\sin(t)} \sin^2(2t) + e^{\cos(t)} \cos^2(2t).$ 

\begin{figure}[ht]
\centering 
\includegraphics[width=0.6\textwidth]{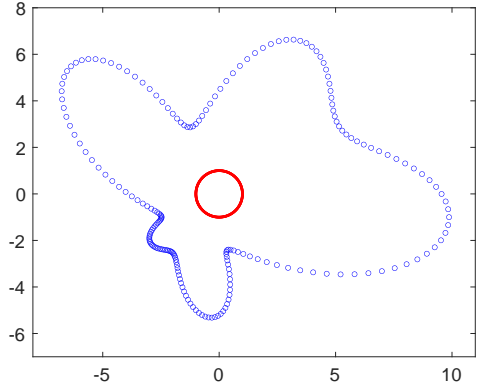}
\caption{ Collocation points on the boundary of the unit disk, marked with $\txtr{\cdot}$ and source points on $\Gamma$ marked with $\txtb{\circ}$.} \label{fig:pontos_teste}
\end{figure}

\begin{figure}[ht]
\centering 
\includegraphics[width=0.48\textwidth]{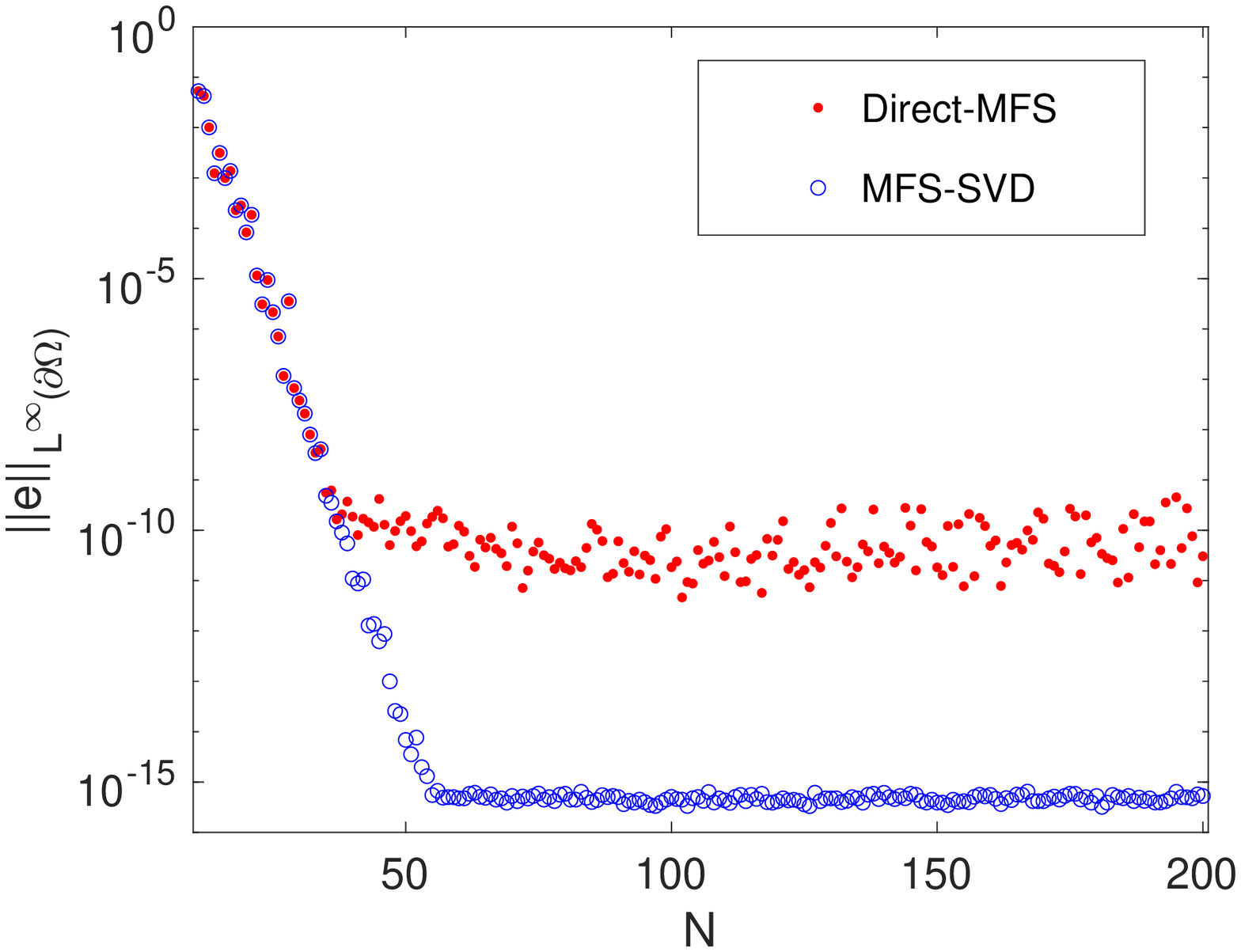}
\includegraphics[width=0.48\textwidth]{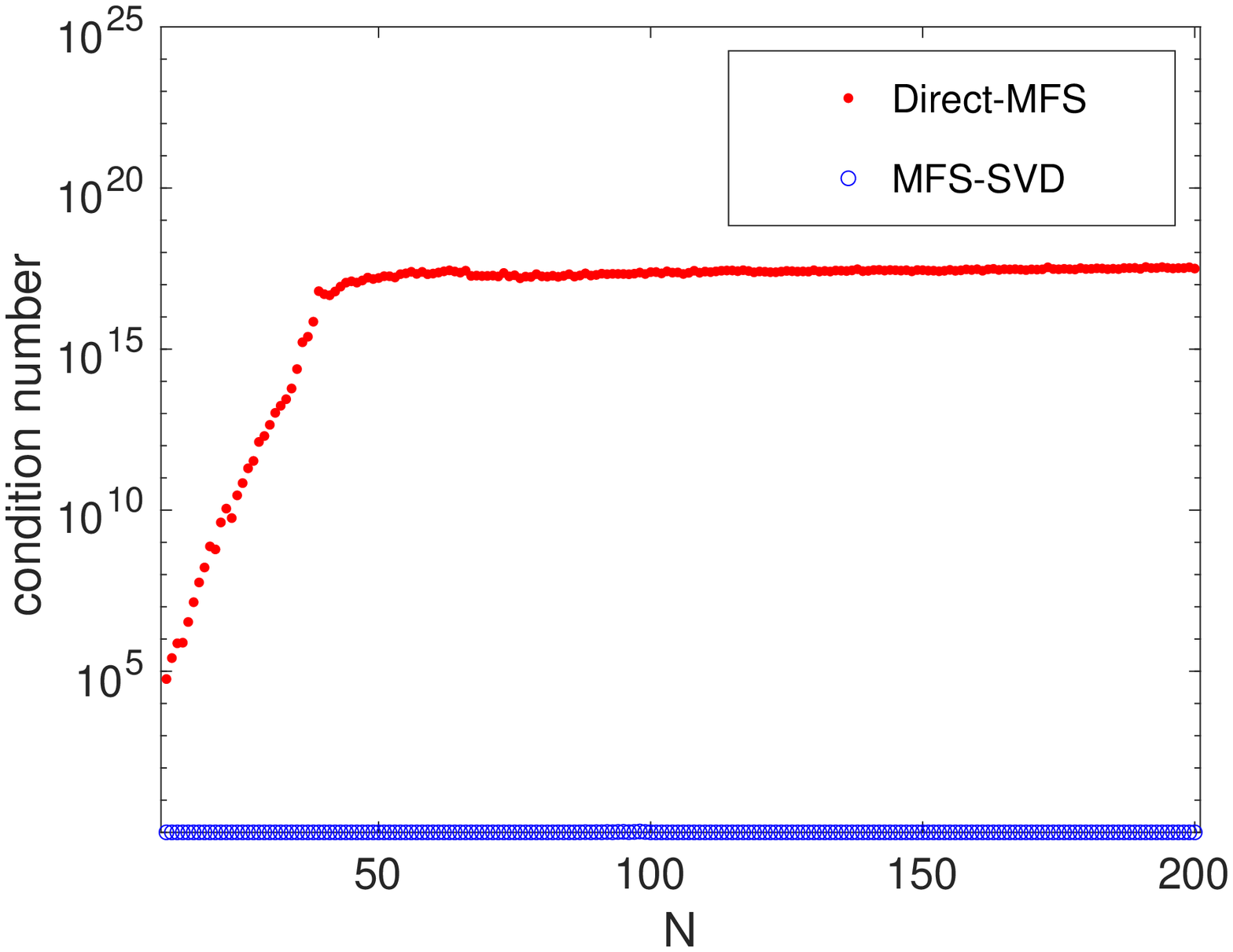}
\caption{Plot of the $L^\infty$ norm of the error on the boundary of the numerical approximations given by the Direct-MFS and MFS-SVD, as a function of $N$ (left plot) and plot of the condition number of the matrix of the linear systems (right plot).} \label{fig:exemplo1}
\end{figure}

Figure~\ref{fig:pontos_teste} shows the collocation points on the boundary of the unit disk, marked with $\txtr{\cdot}$ and the source points on $\Gamma$ marked with $\txtb{\circ}$. Figure~\ref{fig:exemplo1}-left shows the $L^\infty$ norm of the error on the boundary of the numerical approximations given by the Direct-MFS and MFS-SVD, as a function of the number of basis functions, $N$. Again, we took the boundary data given by $g(x,y)=x^2y^3.$ The right plot of the same figure shows the condition number of the matrix of the linear system. We can observe that for $N\leq N_0\approx30$ we obtain similar errors with the Direct-MFS and the MFS-SVD because the basis functions of both approaches span the same functional space. For $N>N_0$ the convergence of the Direct-MFS breaks down due to ill-conditioning, while the errors of the MFS-SVD decrease until we reach values close to machine precision, keeping the condition number to order 1, independently of $N$.

\subsection{An Arnoldi-SVD Method of Fundamental Solutions}
In this section we introduce a technique to avoid the increasing powers of $\frac{r}{R_\Omega}$ in $\mathbf{F}(r,\theta)$, through expansion \eqref{expbc}, for suitable matrices $\mathbf{K}$ and vector valued functions $\mathbf{J}(r,\theta)$. The coefficients of the MFS-SVD are calculated by solving the linear system \eqref{sistema-novo} involving the matrix
\[\left[\boldsymbol{\phi}(r_i,\theta_i)\right]^T=\left[\mathbf{H}(r_i,\theta_i)\right]^T.\left(\mathbf{V}_1^\ast\right)^T.\]
Taking into account \eqref{expbc} we have
\[\left[\mathbf{F}(r_i,\theta_i)\right]^T=\left[\mathbf{H}(r_i,\theta_i)\right]^T\cdot\mathbf{K}^T\]
and defining
\[z_i:=\frac{r_i}{R_\Omega}e^{i\theta_i}\quad\text{and}\quad w_i:=\frac{r_i}{R_\Omega}e^{-i\theta_i}\]
we have
\[\left[\mathbf{F}(r_i,\theta_i)\right]^T=\begin{bmatrix}1 & z_1 & z_1^2 &\cdots & z_1^p & w_1 & w_1^2 & \cdots & w_1^p\\
1 & z_2 & z_2^2 &\cdots & z_2^p & w_2 & w_2^2 & \cdots & w_2^p\\
\vdots & \vdots &\vdots &\ddots & \vdots & \vdots & \vdots & \ddots & \vdots\\
1 & z_M & z_M^2 &\cdots & z_M^p & w_M & w_M^2 & \cdots & w_M^p
\end{bmatrix}:=\left[\begin{array}{c|c}
\mathbf{Z}&\mathbf{W}
\end{array}\right],\]
where the block $\mathbf{Z}$ is a Vandermonde matrix. The matrix $\mathbf{W}$ can also be obtained from a Vandermonde matrix $\tilde{\mathbf{W}}$ after excluding the first column. As pointed out in~\cite{trefethen}, the columns of a Vandermonde matrix, such as $\mathbf{Z},$ can be regarded as vectors $q_0,$ $\mathbf{X}_{\mathbf{Z}} q_0,$ $\mathbf{X}_{\mathbf{Z}}^2 q_0,$ ..., where $q_0=[1\cdots 1]^T$ and $\mathbf{X}_{\mathbf{Z}}=\text{diag}(z_1,\cdots, z_M).$ The Arnoldi process allows us to orthogonalize at each step leading to a sequence of orthogonal vectors $q_0,$ $q_1,$ $q_2,$ ... spanning the same space as the columns of $\mathbf{Z}$. After $n$ steps, we obtain $n+1$ orthogonal vectors $q_0,...,q_n$ and an $(n+1)\times n$ Hessenberg matrix $\mathbf{H}_{\mathbf{Z}}$ that verifies
\[\mathbf{X}_{\mathbf{Z}}\cdot \mathbf{Q}_{-,\mathbf{Z}}=\mathbf{Q}_{\mathbf{Z}}\cdot \mathbf{H}_{\mathbf{Z}},\]
where $\mathbf{Q}_{\mathbf{Z}}$ is a matrix whose columns are the orthogonal vectors $q_0,...,q_n$ and $\mathbf{Q}_{-,\mathbf{Z}}$ is obtained from $\mathbf{Q}_{\mathbf{Z}}$ after removing the last column. The matrices $\mathbf{Z}$ and $\mathbf{Q}_{\mathbf{Z}}$ are related through an upper triangular matrix $\mathbf{R}_{\mathbf{Z}},$ by
\[\mathbf{Z}=\mathbf{Q}_{\mathbf{Z}}\cdot\mathbf{R}_{\mathbf{Z}}.\]
In a similar fashion, for the matrix $\tilde{\mathbf{W}}$, we obtain
\[\tilde{\mathbf{W}}=\mathbf{Q}_{\mathbf{W}}\cdot\tilde{\mathbf{R}}_{\mathbf{W}}\]
and
\[\mathbf{W}=\mathbf{Q}_{\mathbf{W}}\cdot\mathbf{R}_{\mathbf{W}},\]
where $\mathbf{R}_{\mathbf{W}}$ is obtained from $\tilde{\mathbf{R}}_{\mathbf{W}}$ after excluding the first column.

Therefore,
\[\left[\mathbf{F}(r_i,\theta_i)\right]^T=\left[\begin{array}{c|c}
\mathbf{Z}&\mathbf{W}
\end{array}\right]=\left[\begin{array}{c|c}
\mathbf{Q}_{\mathbf{Z}}\cdot\mathbf{R}_{\mathbf{Z}}&\mathbf{Q}_{\mathbf{W}}\cdot\mathbf{R}_{\mathbf{W}}
\end{array}\right]=\left[\begin{array}{c|c}
\mathbf{Q}_{\mathbf{Z}}&\mathbf{Q}_{\mathbf{W}}
\end{array}\right]\cdot\left[\begin{array}{cc}
\mathbf{R}_{\mathbf{Z}}&\mathbf{0} \\
\mathbf{0} & \mathbf{R}_{\mathbf{W}}
\end{array}
\right]\]
and this last equality can be regarded as the evaluation of \eqref{expbc} at the boundary points $(r_i,\theta_i),$ 
\[\mathbf{F}(r_i,\theta_i)=\underbrace{\left[\begin{array}{cc}
\mathbf{R}_{\mathbf{Z}}&\mathbf{0} \\
\mathbf{0} & \mathbf{R}_{\mathbf{W}}
\end{array}
\right]^T}_{=\mathbf{K}}\cdot\underbrace{\left[\begin{array}{c}
\mathbf{Q}_{\mathbf{Z}}^T\\
\mathbf{Q}_{\mathbf{W}}^T
\end{array}
\right]}_{=\mathbf{J}(r_i,\theta_i)}.\]
The evaluation of $\mathbf{J}(r,\theta)$ at a general point $(r,\theta)$ can be performed using the Hessenberg matrices $\mathbf{H}_{\mathbf{Z}}$ and $\mathbf{H}_{\mathbf{W}}$, as in the routine \texttt{polyvalA} of reference \cite{trefethen} and the MFS-SVD solution is evaluated using \eqref{expdisc3}.

The whole numerical technique that we developed in this work is summarized in Algorithm~\ref{numproc}.
\begin{algorithm}
\caption{MFS-SVD}\label{numproc}
\begin{algorithmic}[1]
\State Choose $N\in\mathbb{N}$, the number of basis functions.
\State Place $M(=2N)$ collocation points on $\partial\Omega.$
\State Choose $N$ source points, satisfying the geometric constraint described in section~\ref{neqtech}.
\State Build the matrices $\mathbf{M},$ $\mathbf{Z}$ and $\mathbf{W}.$
\State Calculate $\mathbf{Q}_{\mathbf{Z}},$ $\mathbf{Q}_{\mathbf{W}},$ $\mathbf{H}_{\mathbf{Z}}$ and $\mathbf{H}_{\mathbf{W}}$ using Arnoldi iteration.
\State Calculate $\mathbf{R}_{\mathbf{Z}}$ and $\mathbf{R}_{\mathbf{W}}.$
\State Define $\mathbf{K}.$
\State Calculate $\mathbf{M}_1.$
\State Calculate the SVD factorization of $\mathbf{M}_1$ to obtain $\mathbf{V}_1^\ast.$
\State Calculate the system matrix of the MFS-SVD.
\State Solve the linear system to calculate the MFS-SVD coefficients.
\State Evaluate the solution at an arbitrary point $(r,\theta).$
\end{algorithmic}
\end{algorithm}

\section{Numerical results}
\label{numres}
In this section we present some numerical results to illustrate the performance of the \emph{MFS-SVD}. 

\begin{figure}[ht]
\centering 
\includegraphics[width=0.55\textwidth]{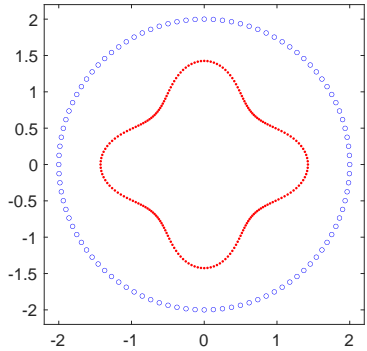}
\caption{Plot of collocation and source points in example 1.} \label{fig:pontos_kara}
\end{figure}

\begin{figure}[ht]
\centering 
\includegraphics[width=0.48\textwidth]{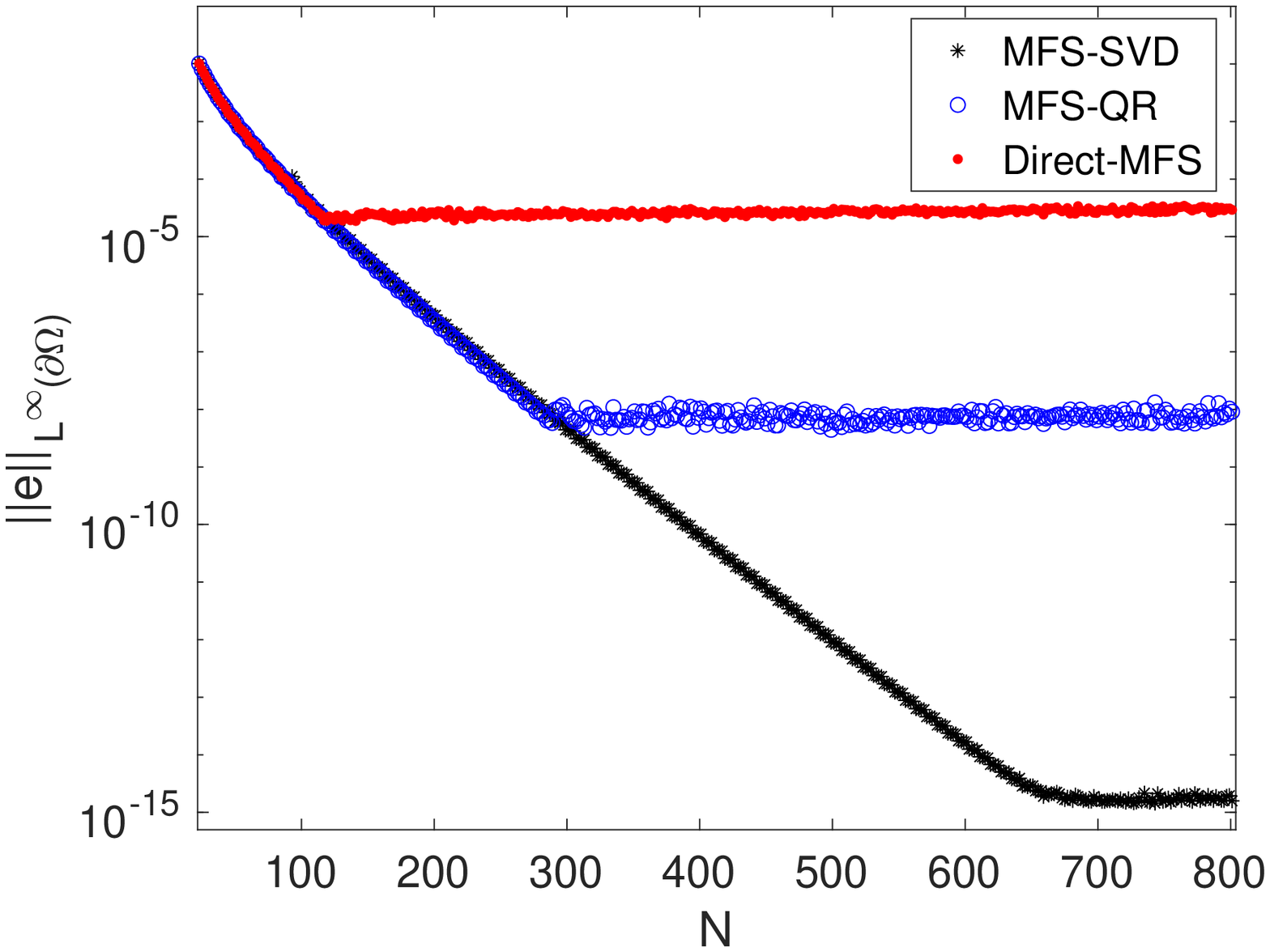}
\includegraphics[width=0.48\textwidth]{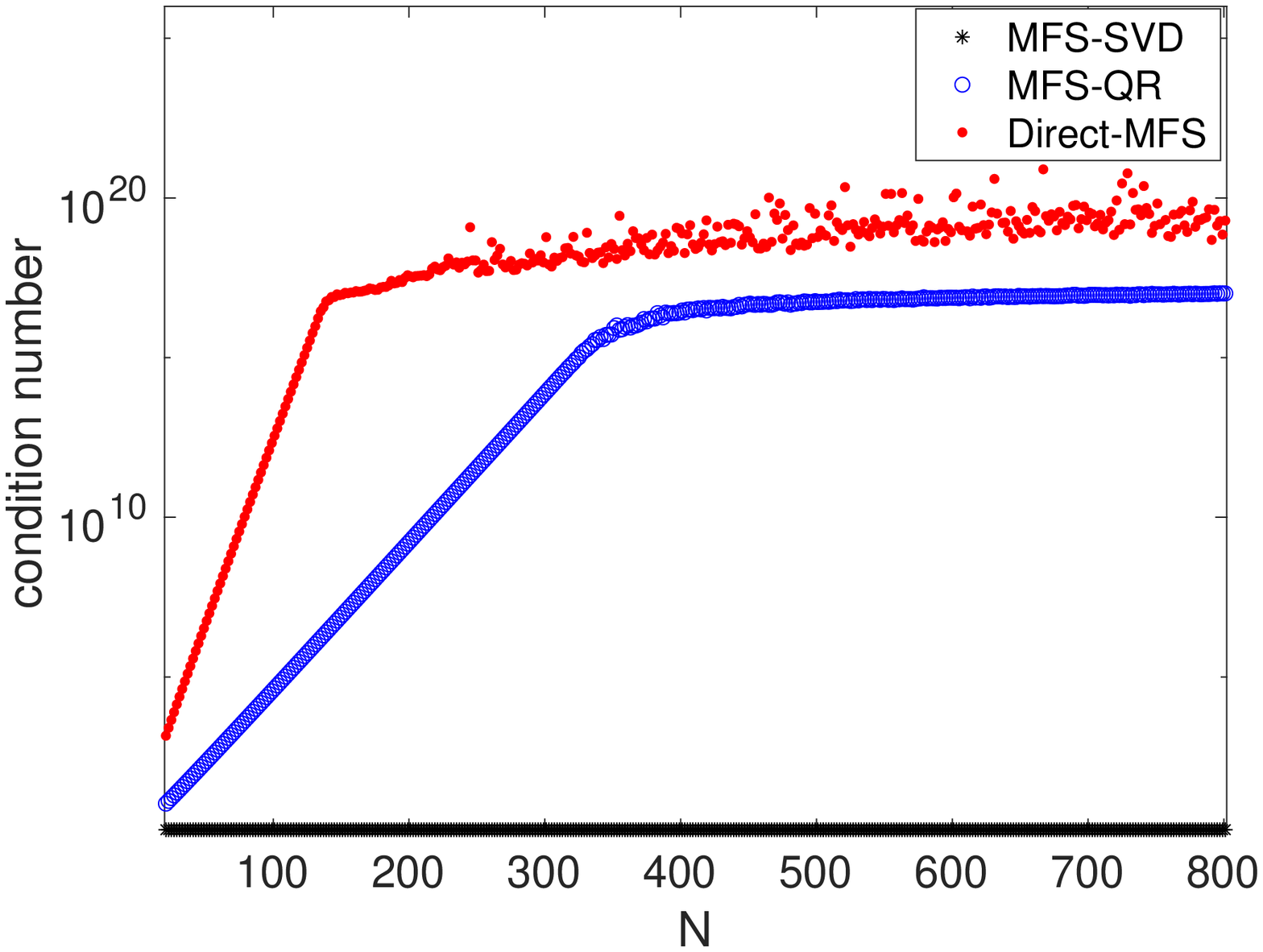}
\caption{Plot of the $L^\infty$ norm of the error of the approximations obtained with Direct-MFS, MFS-QR and MFS-SVD, as a function of $N$ (left plot) and plot of the condition number of the system matrices of the three methods.} \label{fig:kara_resultados}
\end{figure}

In the first numerical example we take the domain with boundary parame\-trized by
\[r(t)=\left(\cos(4t)+\sqrt{\frac{18}{5}-\sin^2(4t)}\right)^{\frac{1}{3}},\ t\in[0,2\pi[\]
also considered in~\cite{Li-CS-Kara} and the same boundary condition defined through the function $g(x,y)=x^2y^3$. In order to allow for a comparison with the MFS-QR we will consider the artificial boundary to be the boundary of a disk centered at the origin and radius equal to 2. Figure~\ref{fig:pontos_kara} shows the collocation points on the boundary and the source points.
Figure~\ref{fig:kara_resultados} shows results obtained with the Direct-MFS, the MFS-QR and the MFS-SVD. The left plot shows the $L^\infty$ norm of the error, measured at 10001 points on the boundary, as a function of $N,$ while the right plot of the same figure shows the condition number of the matrices of the linear systems. We can observe that the numerical results obtained with the three approaches are similar for $N\leq N_0\approx120$ because the functional spaces are the same although defined through different basis functions. For $N>N_0,$ the convergence of the Direct-MFS breaks down due to ill-conditioning. Note that the condition numbers corresponding to the Direct-MFS grow exponentially. The condition numbers corresponding to the MFS-QR also grow exponentially, but at a lower rate than the Direct-MFS. This growth is due to the increasing powers of the harmonic polynomials in the MFS-QR expansion. The numerical results of the MFS-QR and MFS-SVD are similar for $N\leq N_1\approx290.$ For $N>N_1$ the convergence of the MFS-QR also breaks down due to ill-conditioning, while the MFS-SVD converges until we reach accuracy close to machine precision, keeping the condition number approximately equal to 1.65, independently of $N$.

\begin{figure}[ht]
\centering 
\includegraphics[width=0.6\textwidth]{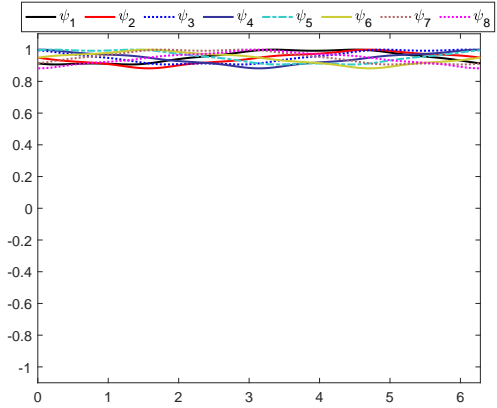}
\caption{Plot of the basis functions of the Direct-MFS, when $N=8$ and $R=10$.} \label{fig:fbmfs}
\end{figure}

\begin{figure}[ht]
\centering 
\includegraphics[width=0.48\textwidth]{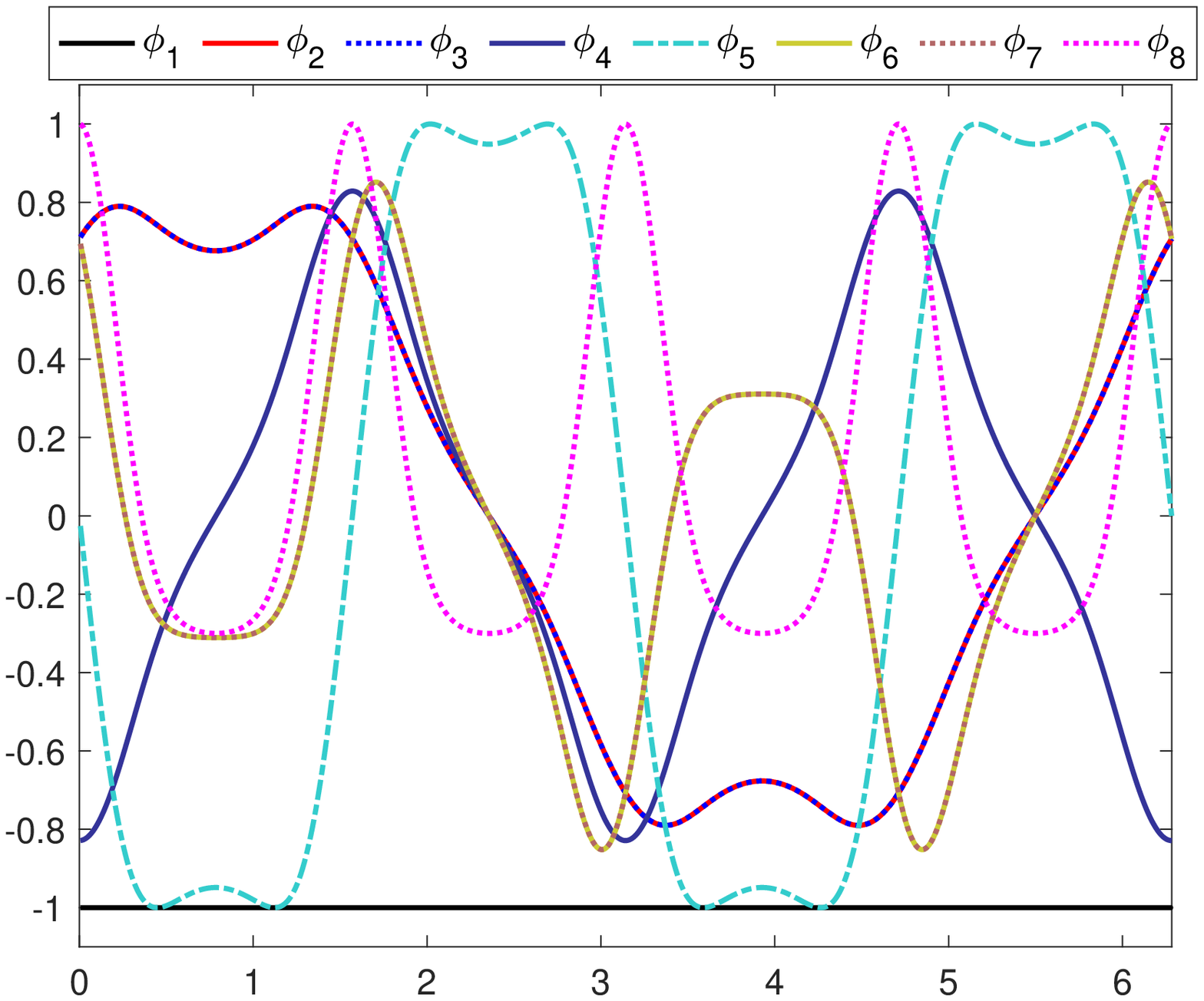}
\includegraphics[width=0.48\textwidth]{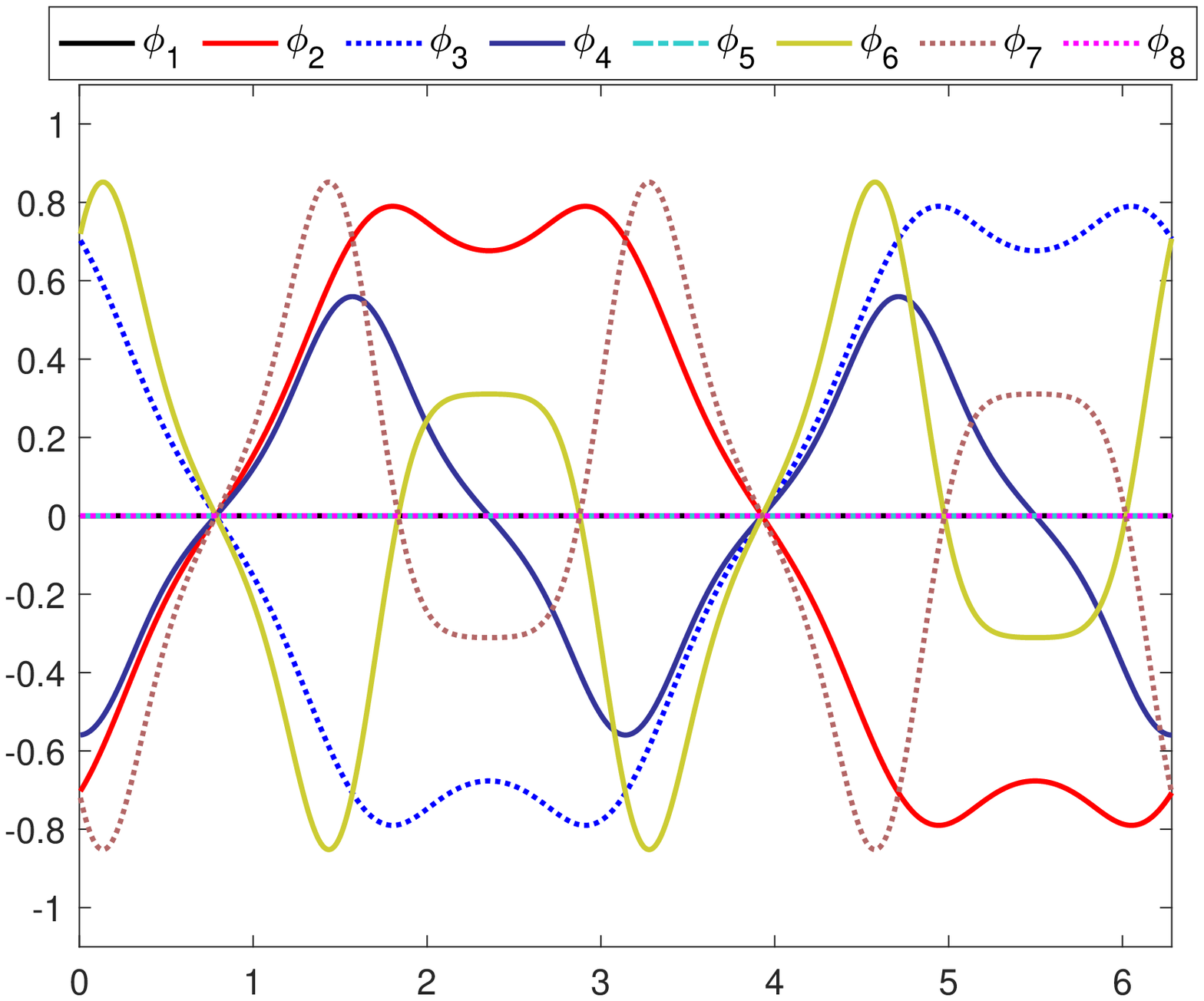}
\caption{Plots of the real part (left plot) and imaginary part (right plot) of the (well conditioned) basis functions corresponding to MFS-SVD spanning the same functional space as the Direct-MFS basis functions plotted in Figure~\ref{fig:fbmfs}.} \label{fig:fbmfssvd}
\end{figure}

The ill-conditioning of the Direct-MFS is related to the fact that all the basis functions become almost linearly dependent on $\partial\Omega$. This effect is even more evident if we take the artificial boundary \textit{far} from $\partial\Omega.$ Figure~\ref{fig:fbmfs} shows the plots of $N=8$ Direct-MFS basis functions, normalized to have unit $L^\infty$ norm, for source points equally spaced on the boundary of the disk with radius equal to 10. All the basis functions become almost undistinguishable at this scale. 
Figure~\ref{fig:fbmfssvd} shows the real part (left plot) and imaginary part (right plot) of the (well conditioned) basis functions corresponding to MFS-SVD spanning the same functional space as the Direct-MFS basis functions.

\begin{figure}[ht]
\centering 
\includegraphics[width=0.48\textwidth]{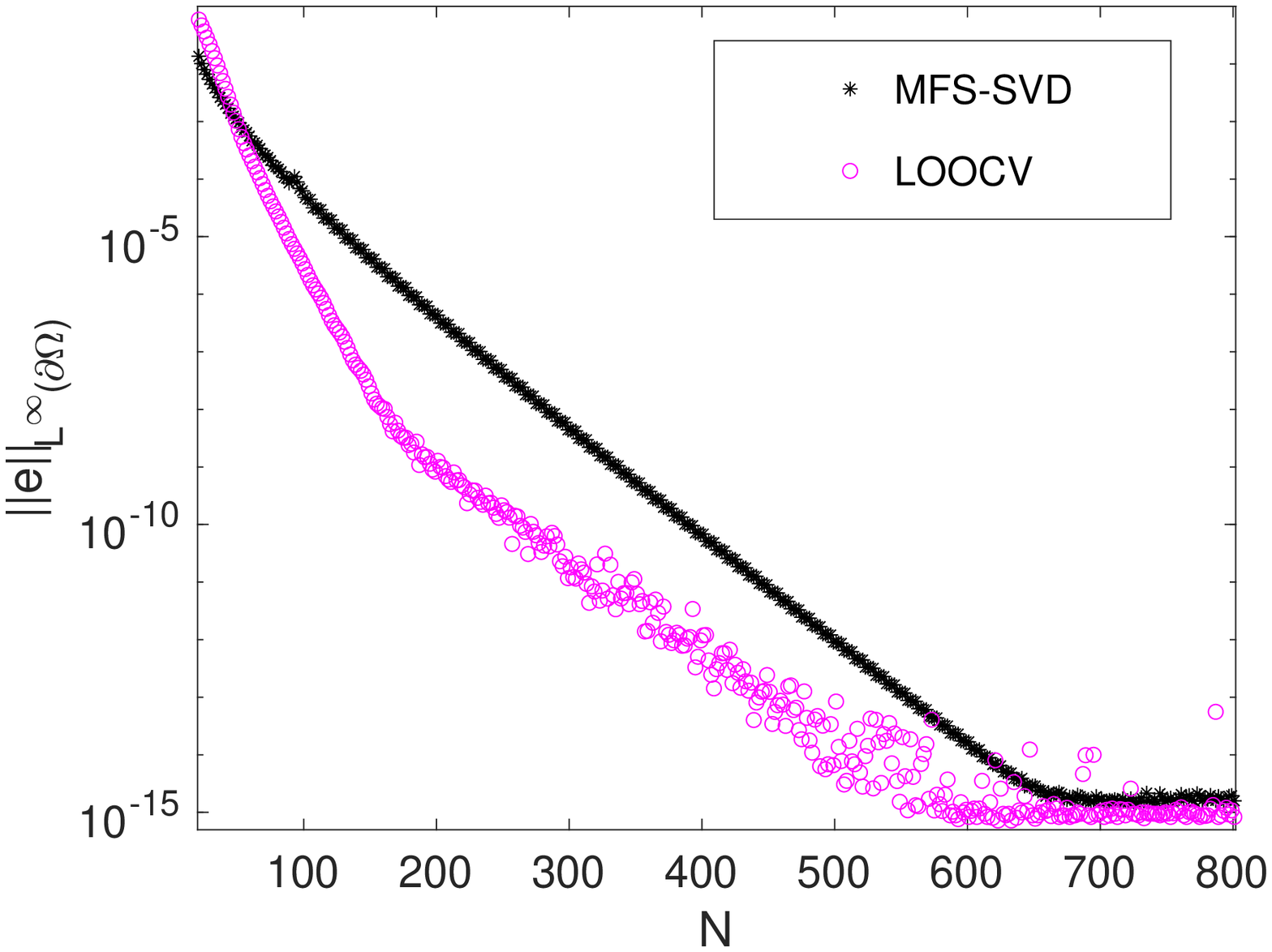}
\caption{Plot of the $L^\infty$ norm of the error of the approximations obtained with the  MFS-SVD and LOOCV, as a function of $N$.} \label{fig:kara_loocv}
\end{figure}

\begin{figure}[ht]
\centering 
\includegraphics[width=0.48\textwidth]{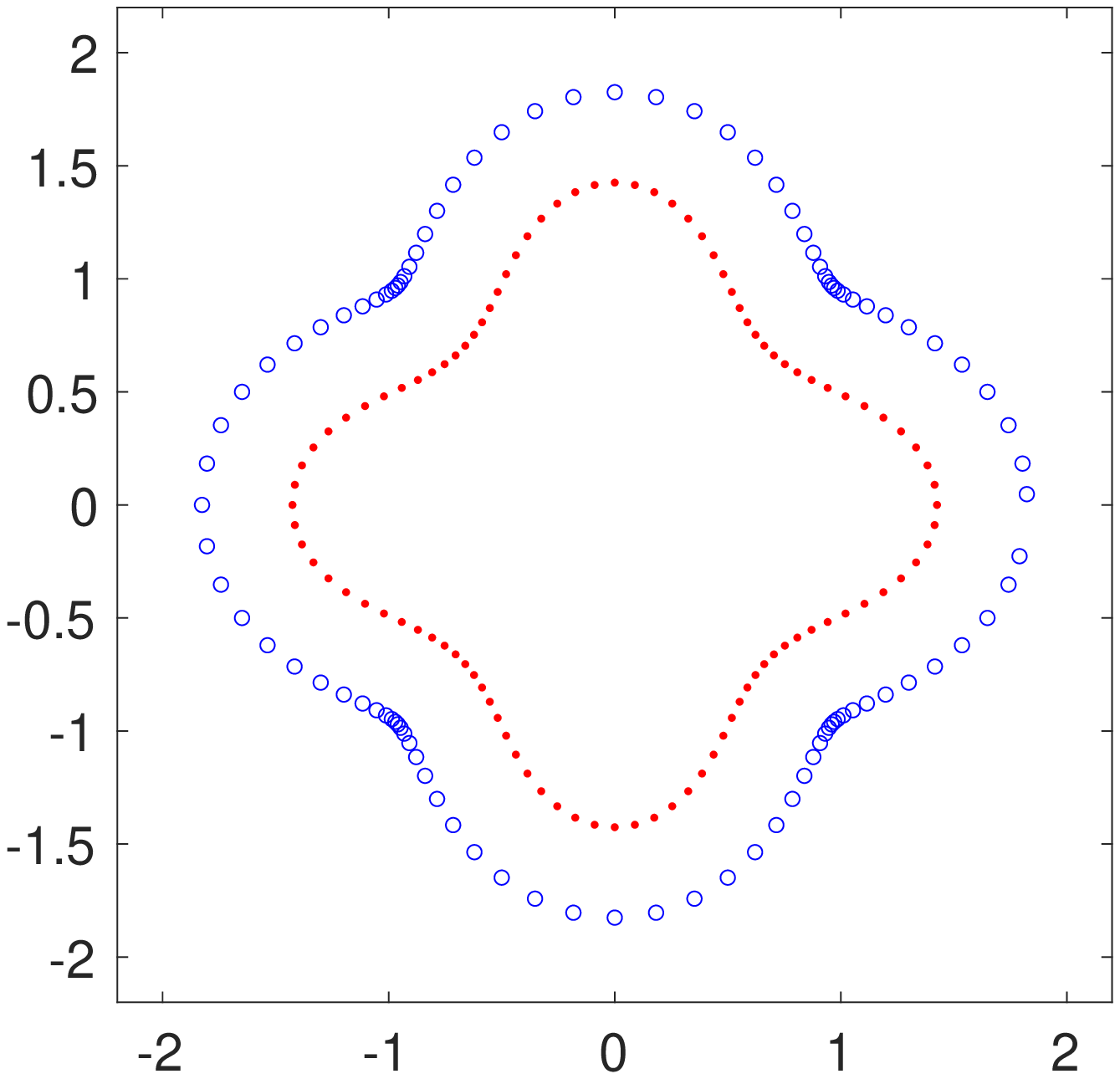}
\includegraphics[width=0.48\textwidth]{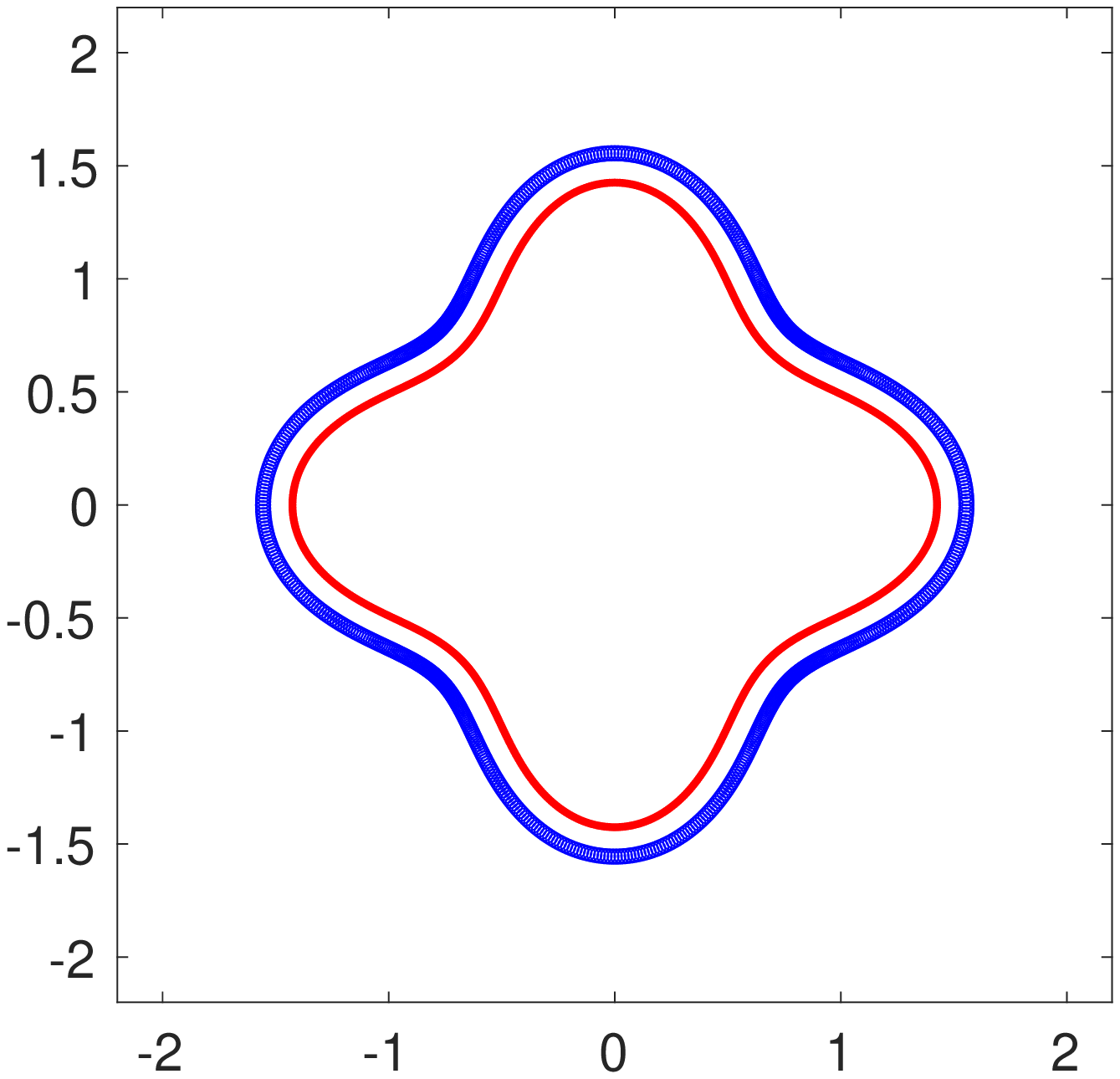}
\caption{Collocation and source points obtained with LOOCV technique, when $N=100$ (left plot) and $N=800$ (right plot).} \label{fig:kara_loocv_pts}
\end{figure}

In the previous example we took the artificial boundary to be a circle, in order to compare with the results obtained with the MFS-QR. However, the numerical results obtained suggest that placing the source points on a circle might be not a convenient choice. An interesting question would be to find the optimal location for placing the source points which should depend on the geometry of the domain and on the boundary data. One possibility is to use an adaptive technique that for fixed domain and boundary data optimize the location of the source points, such as the LOOCV algorithm proposed by Rippa in~\cite{Rippa} and used in~\cite{Chen-Kara} in the context of the application of the MFS. Figure~\ref{fig:kara_loocv} shows the convergence curve obtained with MFS-SVD and LOOCV techniques. We can observe that the convergence of the LOOCV is faster than the MFS-SVD. Moreover, the results obtained suggest that when the number of collocation points become large, the sources should be getting closer to the boundary. For instance, in Figure~\ref{fig:kara_loocv_pts} we show the collocation and source points obtained with LOOCV technique, when $N=100$ (left plot) and $N=800$ (right plot).

\begin{figure}[ht]
\centering 
\includegraphics[width=0.6\textwidth]{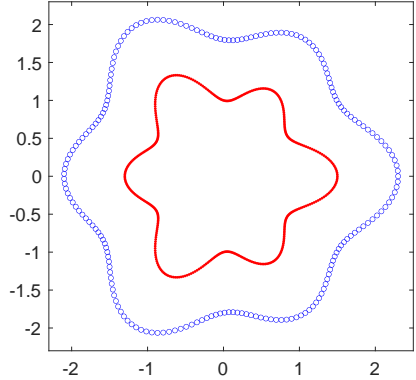}
\caption{Collocation points on the boundary of the unit disk, marked with $\txtr{\cdot}$ and source points marked with $\txtb{\circ}$.} \label{fig:pontos_oscilante}
\end{figure}

In the second numerical example we take the domain with boundary defined by
\[\left\{r_1(t)(\cos(t),\sin(t)),\ t\in[0,2\pi[\right\},\]
where
\[r_1(t)=\frac{6}{5}+\frac{\cos(6t)}{5}+\frac{\cos(3t)}{10},\]
the artificial boundary parametrized by
\[r(t)=2+\frac{\cos(6t)}{5}+\frac{\cos(3t)}{10},\ t\in[0,2\pi[\]
and the same boundary condition defined by the function $g(x,y)=x^2y^3.$ 
Figure~\ref{fig:pontos_oscilante} shows the collocation and source points in this second numerical example.

\begin{figure}[ht]
\centering 
\includegraphics[width=0.48\textwidth]{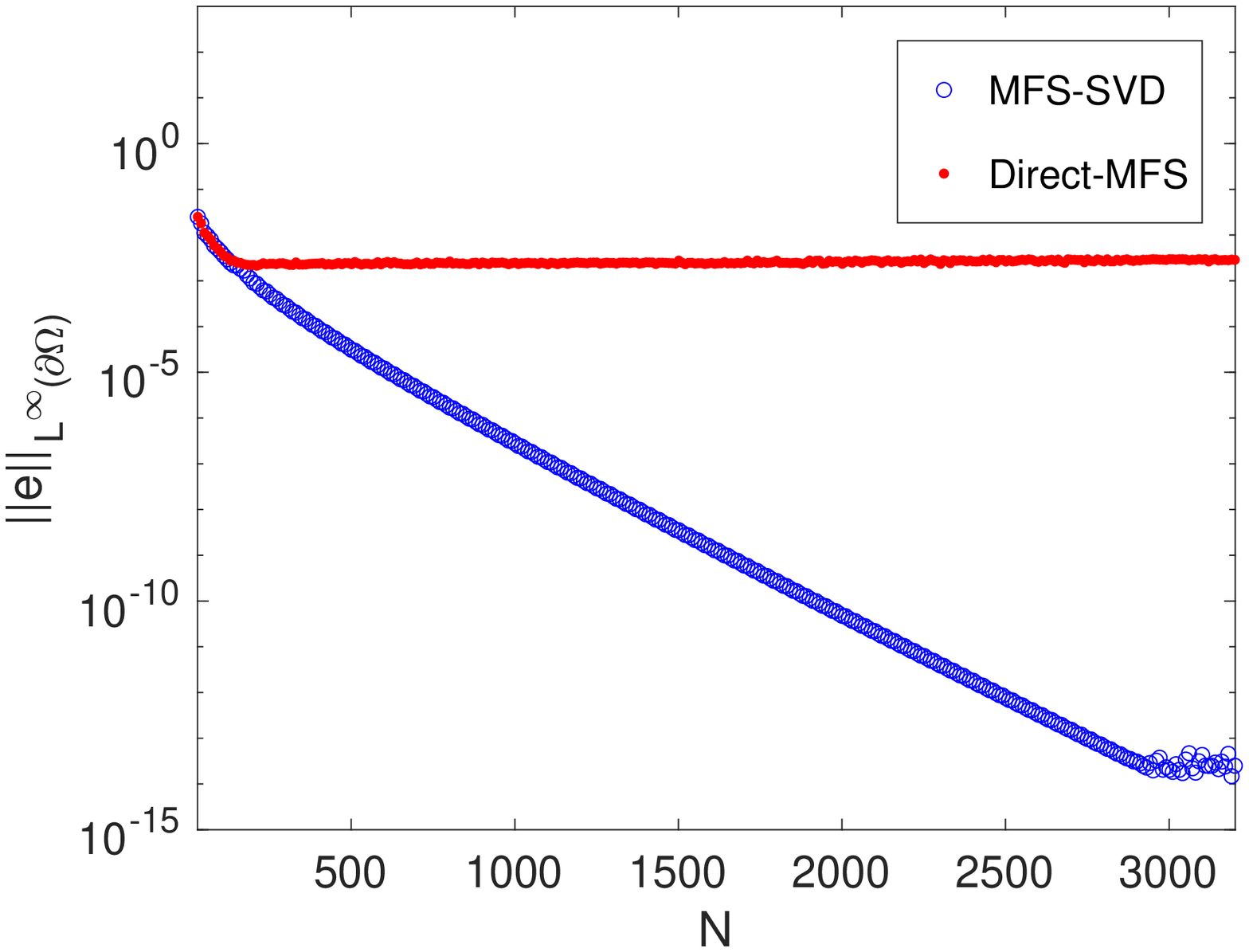}
\includegraphics[width=0.48\textwidth]{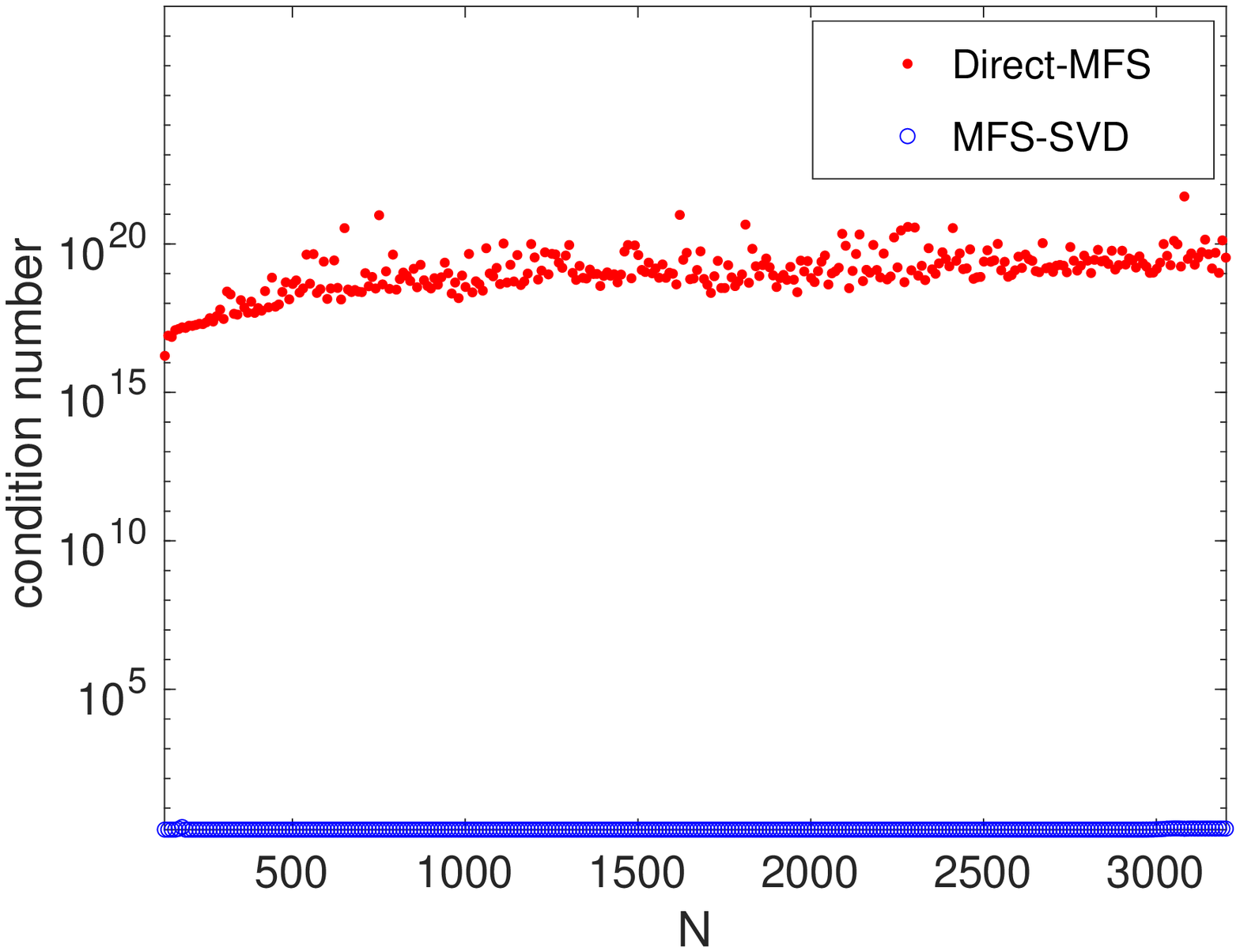}
\caption{Plot of the $L^\infty$ norm of the error of the approximations obtained with Direct-MFS and MFS-SVD, as a function of $N$ in the second numerical example (left plot) and plot of the condition number of the matrices of these two approaches. The boundary condition is defined by the function $g(x,y)=x^2y^3.$} \label{fig:ex2_resultadosprim}
\end{figure}

Figure~\ref{fig:ex2_resultadosprim} shows the plots of the $L^\infty$ norms of the error and condition numbers, both as a function of $N.$ Again, we can observe that the Direct-MFS converges, as we increase $N$ until some $N_0\approx150.$ Then, the convergence breaks down due to ill-conditioning, while the MFS-SVD converges until an accuracy close to machine precision. The condition number is approximately equal to 1.89, independently of $N.$

\begin{figure}[ht]
\centering 
\includegraphics[width=0.6\textwidth]{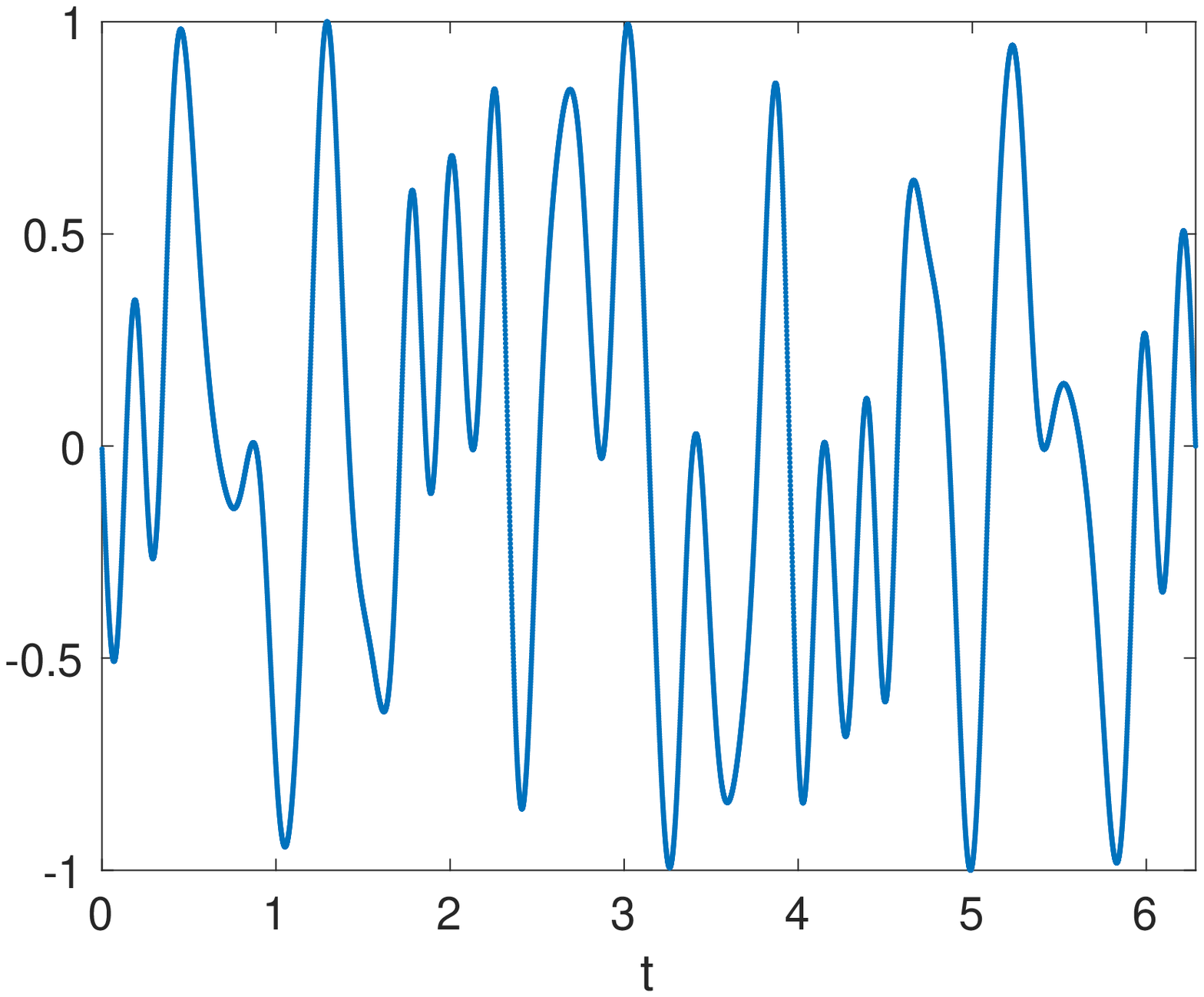}
\caption{ Plot of the boundary data $g(r_1(t)\cos(t), r_1(t)\sin(t)),\ t\in[0,2\pi[.$} \label{fig:bc_oscilante}
\end{figure}

\begin{figure}[ht]
\centering 
\includegraphics[width=0.6\textwidth]{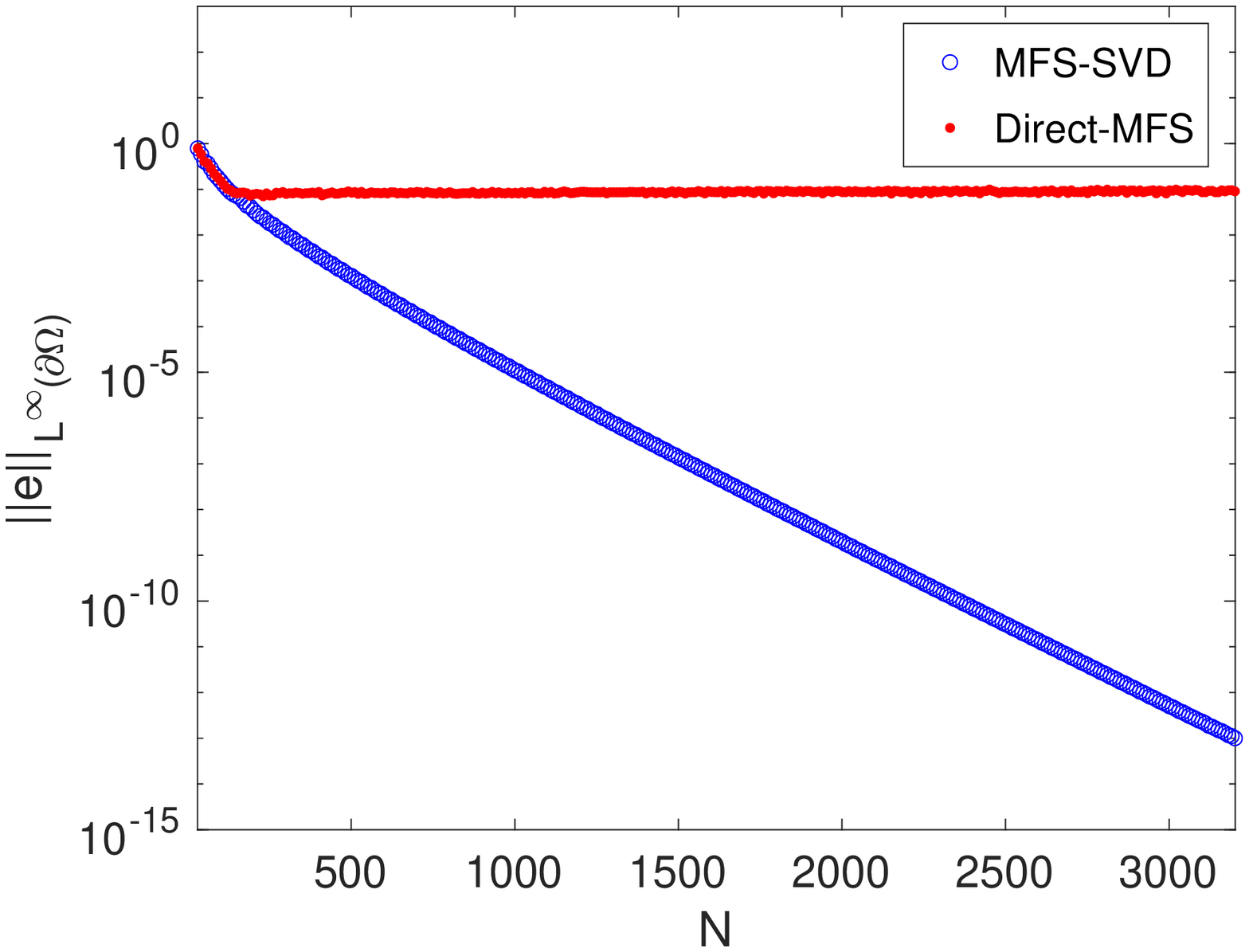}
\caption{Plot of the $L^\infty$ norm of the error of the approximations obtained with Direct-MFS and MFS-SVD, as a function of $N$ in the second numerical example for $g(x,y)=\cos(10x)\sin(10y)$.} \label{fig:ex2_resultados}
\end{figure}

\begin{figure}[ht]
\centering 
\includegraphics[width=0.7\textwidth]{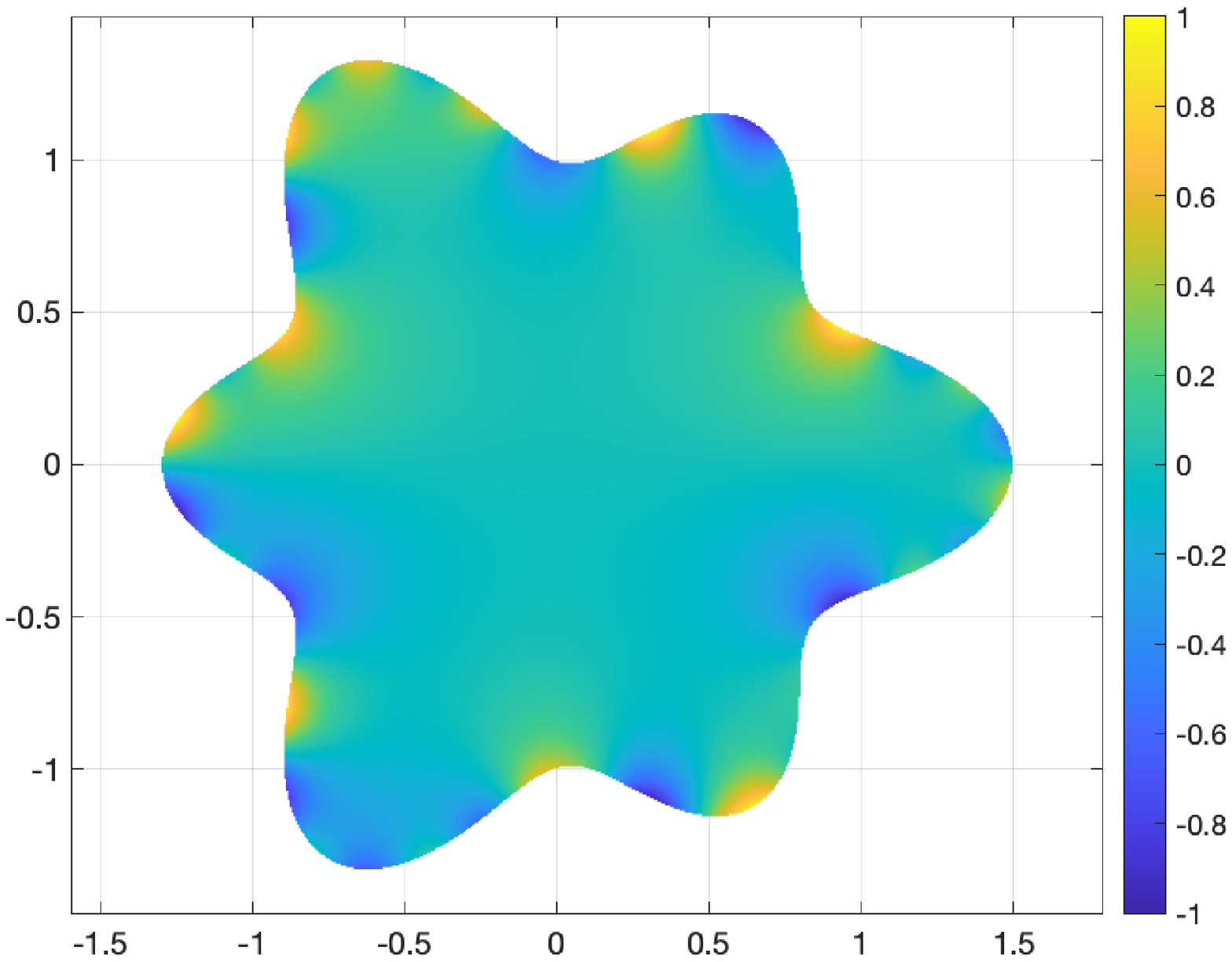}
\caption{Plot of the solution of the boundary value problem of the second numerical example with boundary condition defined by the function $g(x,y)=\cos(10x)\sin(10y)$.} \label{fig:solucao}
\end{figure}

Next, we change the boundary condition to a more oscillatory function
\[g(x,y)=\cos(10x)\sin(10y).\] 
In Figure~\ref{fig:bc_oscilante} we plot the boundary data $g(r_1(t)\cos(t), r_1(t)\sin(t)),\ t\in[0,2\pi[.$
Figure~\ref{fig:ex2_resultados} shows the plots of the $L^\infty$ norm of the error of the approximations obtained with Direct-MFS and MFS-SVD, as a function of $N$ (left plot).
Figure~\ref{fig:solucao} shows the plot of the solution of the boundary value problem of the second numerical example.

\begin{figure}[ht]
\centering 
\includegraphics[width=0.48\textwidth]{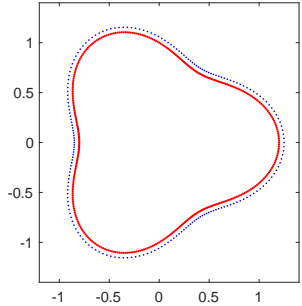}
\caption{The location of the collocation points on the boundary (marked in red) and the source points obtained with $\rho=0.05$ (marked in blue).} \label{fig:exeptos}
\end{figure}

\begin{figure}[ht]
\centering 
\includegraphics[width=0.48\textwidth]{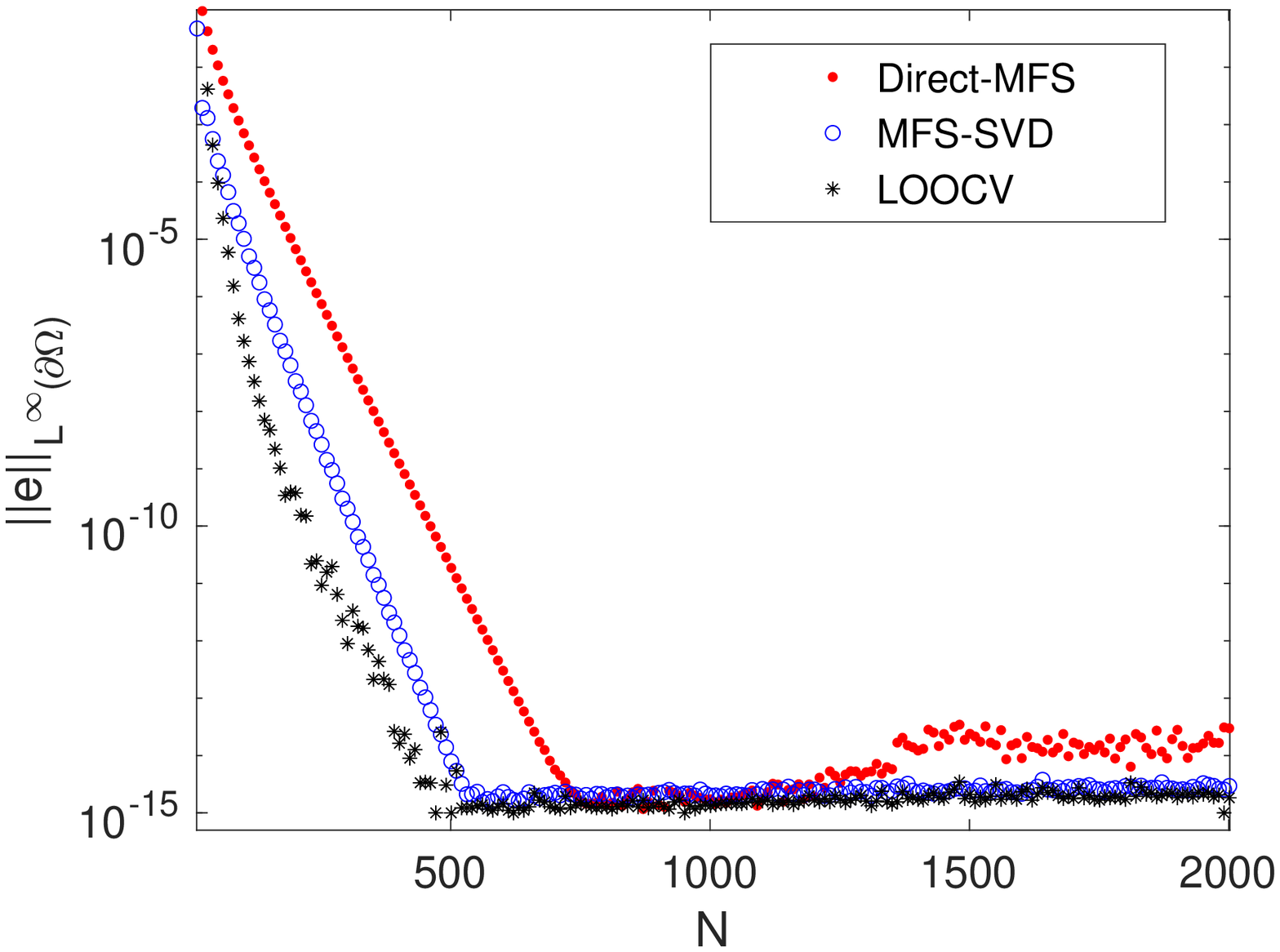}
\includegraphics[width=0.48\textwidth]{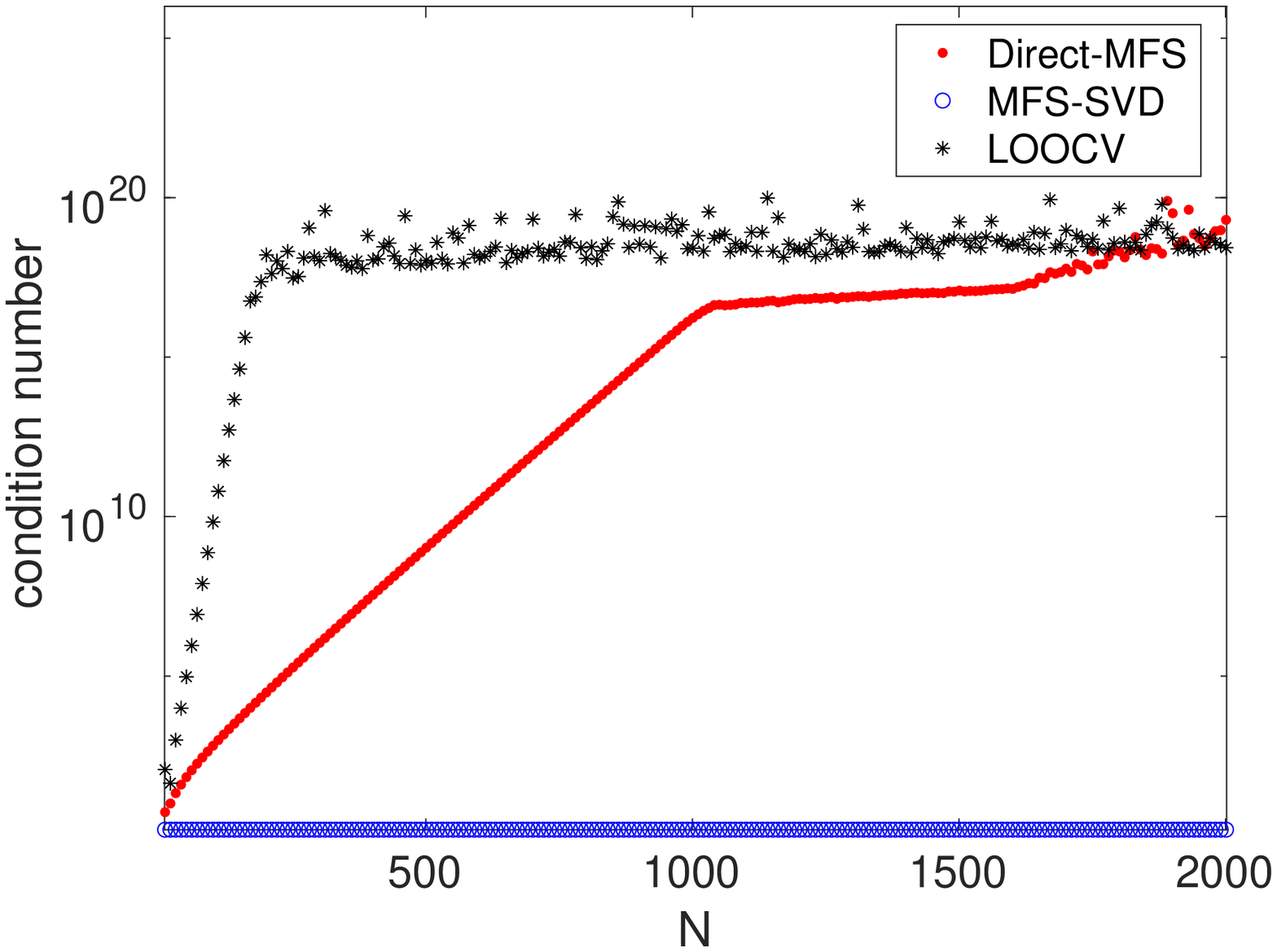}

\caption{Plot of the $L^\infty$ norm of the error of the approximations obtained with Direct-MFS, MFS-SVD and LOOCV, as a function of $N$ (left plot) and plot of the condition number (right plot). The boundary condition is defined by the function $g(x,y)=x^2y^3.$} \label{fig:novex}
\end{figure}

\begin{figure}[ht]
\centering 
\includegraphics[width=0.48\textwidth]{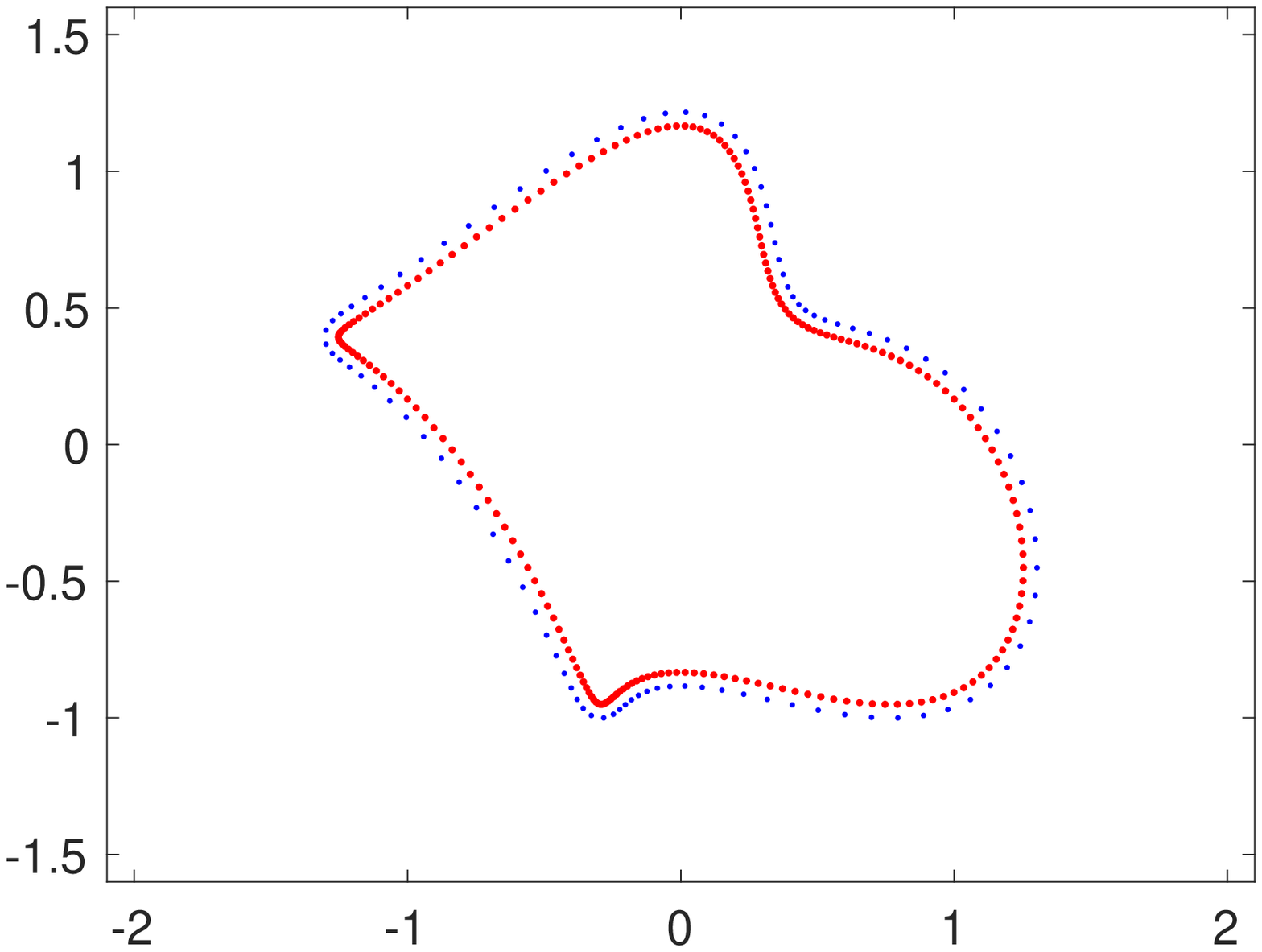}
\includegraphics[width=0.48\textwidth]{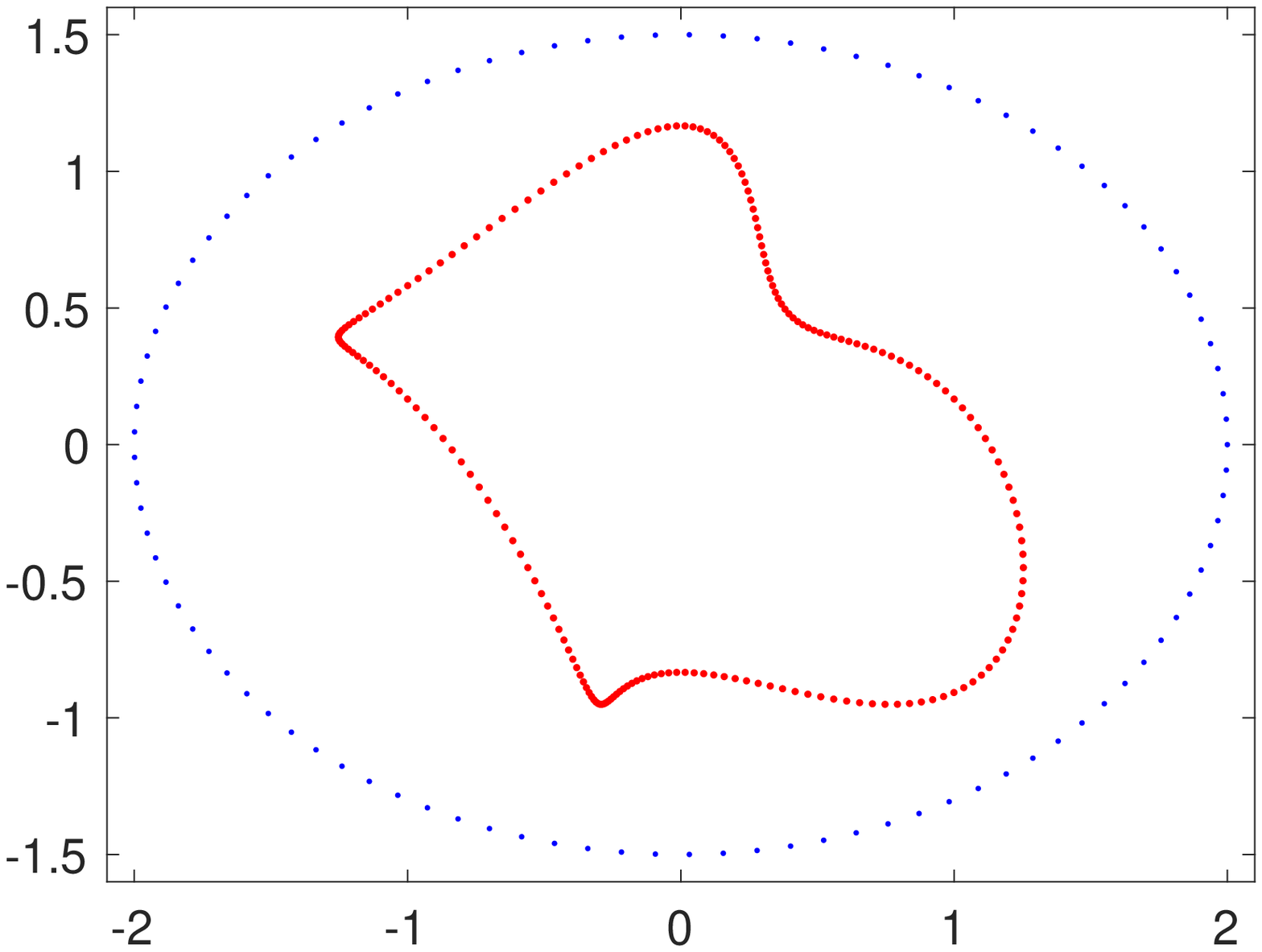}
\caption{The location of the collocation points on the boundary (marked in red) and the source points obtained with $\rho=0.05$ (marked in blue).} \label{fig:pontos_ulti}
\end{figure}

\begin{figure}[ht]
\centering 
\includegraphics[width=0.48\textwidth]{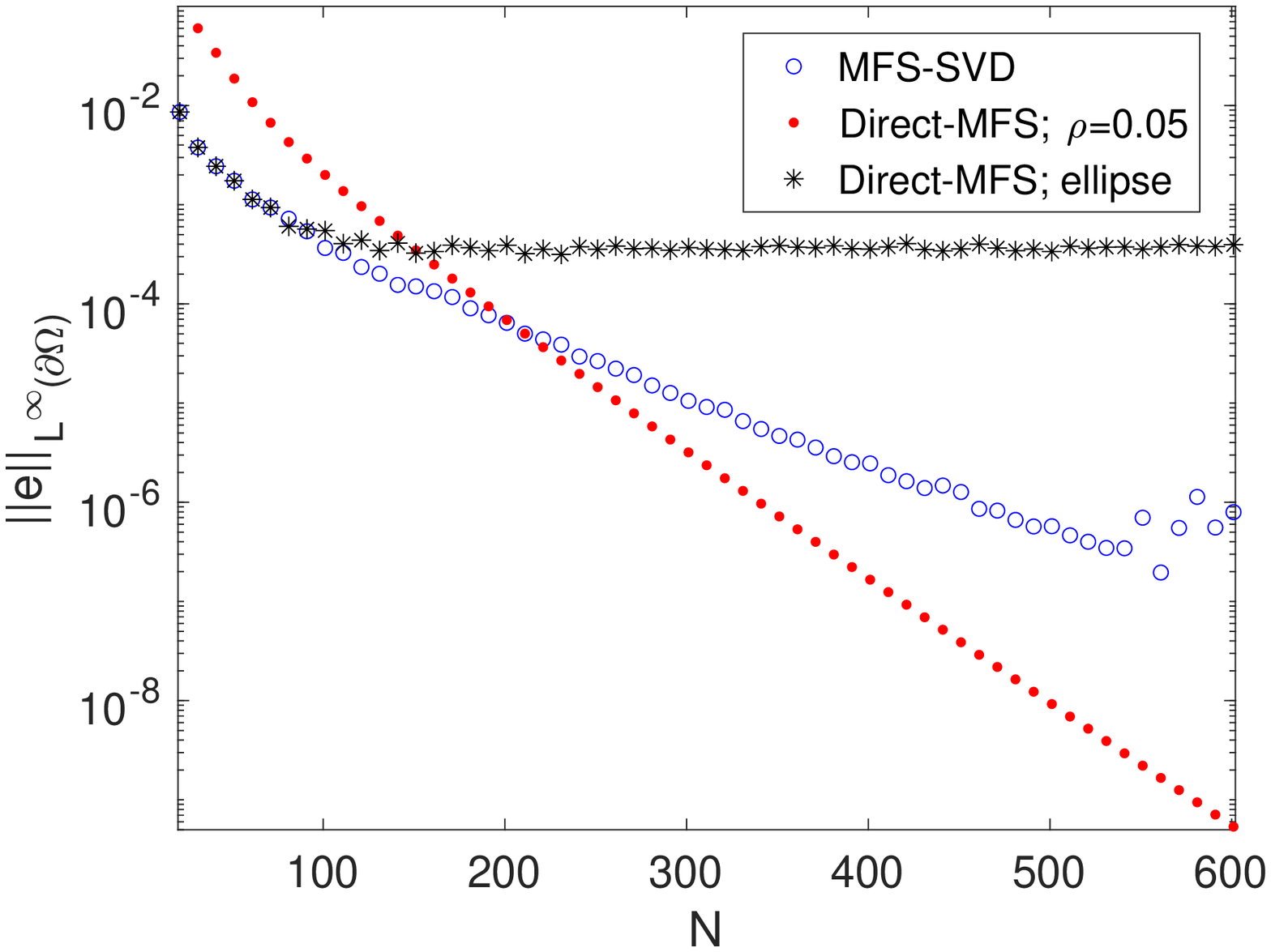}
\includegraphics[width=0.48\textwidth]{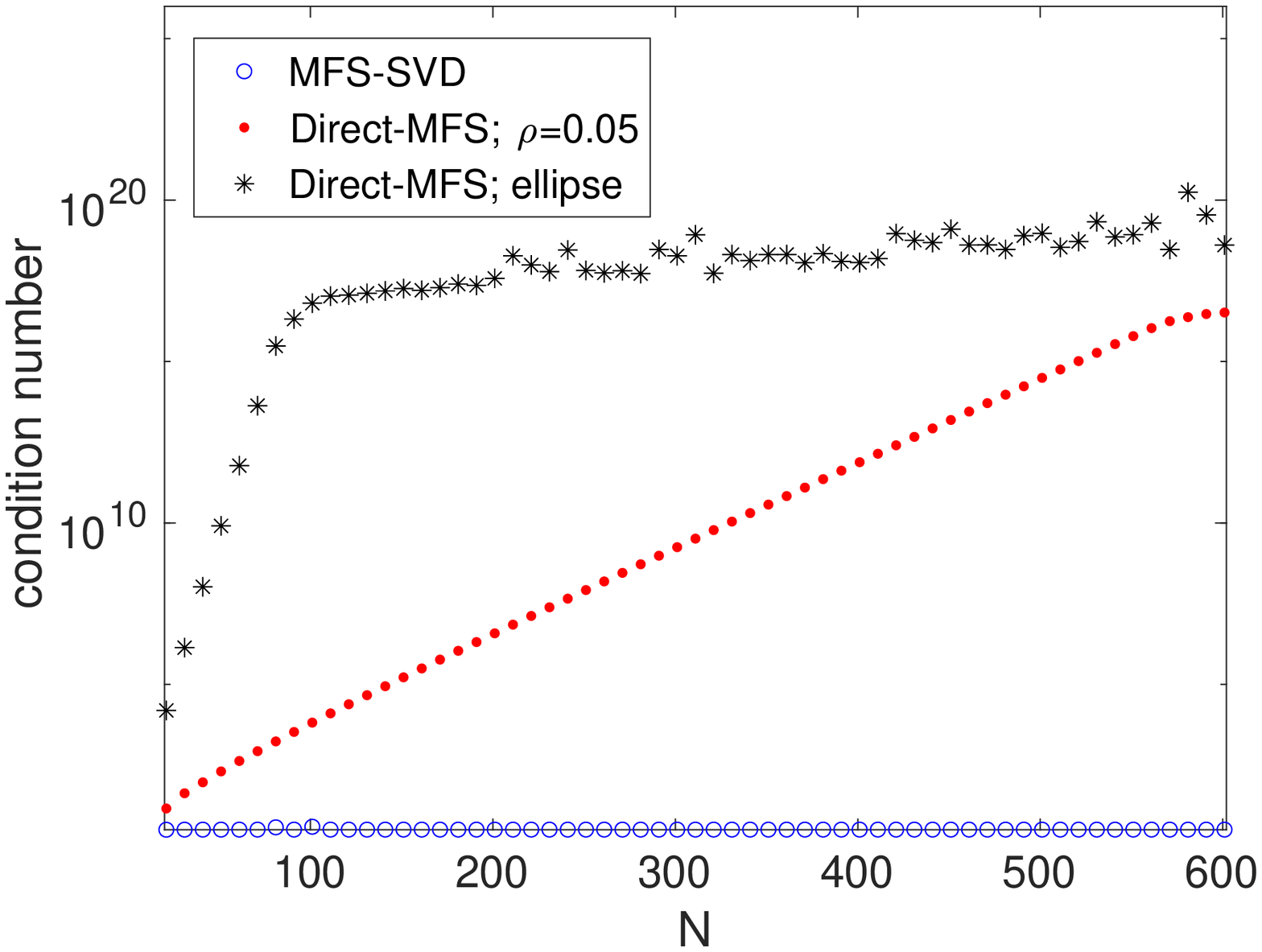}
\caption{Plot of the $L^\infty$ norm of the error of the approximations obtained with Direct-MFS and MFS-SVD, as a function of $N$ in the second numerical example (left plot) and plot of the condition number of the matrices of these two approaches. The boundary condition is defined by the function $g(x,y)=x^2y^3.$} \label{fig:ex2_resultadoulti}
\end{figure}

As reported in literature, eg.~\cite{Chen-Kara}, when solving Laplace equation with non-harmonic boundary conditions, the source points in the Direct-MFS shall be placed close to the physical boundary. However, this $\emph{optimality}$ may be related with the fact that placing the source-points close to the boundary decreases significantly the condition number of the system matrices and this may justify the high accuracy obtained in that case. Placing the source-points far from the physical boundary will always lead to highly ill-conditioned matrices which may prevent the method from achieving high accuracy. On the other hand, the MFS-SVD, as introduced in this work assumes the geometric constraint defined in section~\ref{neqtech} and, in general, this constraint excludes the possibility of the choosing the artificial boundary very close to the physical boundary. The last two numerical examples show performance of the Direct-MFS and the MFS-SVD, but with different artificial boundaries in the two cases. Next, we consider the domains with boundaries defined by
\[\eta_1(t)=\left\{\left(1+\frac{1}{5}\cos(3t)\right)(\cos(t),\sin(t)),\ t\in[0,2\pi[\right\}\]
and
\[\eta_2(t)=\left\{(\cos(t)-\frac{\cos(t)\sin(2t)}{2},\sin(t)+\frac{\cos(4t)}{6}),\ t\in[0,2\pi[\right\}\]
and placing the source-points on the artificial boundary defined by
\begin{equation}
\label{pf}
\tilde{\eta}_i(t)=\eta_i(t)+\rho n_i(t),\ i=1,2
\end{equation}
where $\rho$ is a small positive parameter and $n(t)$ is the unitary outward vector normal to the boundary at the point $\eta_i(t),\ i=1,2$. Figures~\ref{fig:exeptos} and~\ref{fig:pontos_ulti}-left illustrate the choice of the collocation points on the boundary (marked in red) and the source points obtained with $\rho=0.05$ (marked in blue) used in the Direct-MFS. We will apply the MFS-SVD placing source points, respectively, on a circle centered at the origin with radius equal to 1.5 and on the artificial boundary
\[\left\{(2\cos(t),1.5\sin(t)):\ t\in[0,2\pi[\right\},\]
as in figure~\ref{fig:pontos_ulti}-right. Figure~\ref{fig:novex} shows convergence results obtained with three different approaches: 
\begin{enumerate}
\item Direct-MFS with the choice of source points plotted in Figure~\ref{fig:exeptos}, 
\item MFS-SVD with source points placed on a circle centered at the origin with radius equal to 1.5 and 
\item LOOCV algorithm.
\end{enumerate} 
The right plot of the same Figure shows the condition number of the matrix involved in each technique. Note that for the Direct-MFS and MFS-SVD this is the matrix of the linear system, while for LOOCV this is a matrix to be inverted in order to calculated the coefficients of the MFS linear combination. We can observe that, in this case, the convergence of MFS-SVD is faster than the Direct-MFS with small condition numbers, independently of $N$, while the condition number of the Direct-MFS grows exponentially. The convergence of the LOOCV is even faster, but the condition number of the matrix to be inverted increases exponentially. Figure~\ref{fig:ex2_resultadoulti} shows results obtained for $\eta_2$ with the Direct-MFS, MFS-SVD and Direct-MFS with the same source points of the MFS-SVD, plotted in figure~\ref{fig:pontos_ulti}-right. In this case, the Direct-MFS with source points defined by \eqref{pf} with $\rho=0.05$ allows to obtain better results than the MFS-SVD.

\section{Conclusions and future work}
We proposed a new algorithm for generating a set of functions spanning the same functional space as the Direct-MFS basis functions, but which are much better conditioned. This approach allows us to remove the ill-conditioning of the Direct-MFS for general star shaped domains. The artificial boundary for the MFS was assumed to satisfy a certain geometric constraint, defined in section~\ref{neqtech} for ensuring convergence of the series expansions considered. In general, this constraint implies that for the source points cannot be chosen too close to the physical boundary. This is a drawback of the MFS-SVD compared with the Direct-MFS. We hope to address the extension to the general case, where we drop the constraint in the MFS-SVD in a future work. The computational cost of the MFS-SVD is dominated by the calculation of a singular value decomposition and Arnoldi algorithm. Clearly, this computational cost is higher than the Direct-MFS. On the other hand, the MFS-SVD allows obtain well-conditioned system matrices, independently of the choice of artificial boundaries, provided they satisfy the geometric constrain defined in section~\ref{neqtech} and this is a clear advantage when compared with the Direct-MFS, for instance, when solving inverse problems, where the boundary data is affected by noise. Another approach is to consider the LOOCV technique that allows an adaptive technique for determining the source points. Typically, this approach may be even more accurate than the MFS-SVD, but implies the explicit calculation of the inverse of an ill conditioned matrix.

As was mentioned in the Introduction, in our opinion there are two major issued in the MFS to be resolved: a clear criteria for placing the source points and the ill-conditioning. Actually both issues are related in the sense that the \emph{optimal} locations for the source points advocated in the literature were proposed in the context of the application of the Direct-MFS, where the ill-conditioning limits the accuracy that can be achieved and imposes a trade-off between accuracy and conditioning. It is our belief that the MFS-SVD or any other technique to perform a change of basis will allow to explore other choices for artificial boundaries without the constraints imposed by the ill-conditioning. In particular, we believe that it would be interesting to revisit the problem of the location of the source points in the context of a better conditioned basis for the MFS for which the ill-conditioning is not an issue. The three-dimensional case and the extensions to other PDEs are under current research.

\section{Conflict of Interests}
 The author declares that he has no conflict of interest.

\section{Data availability statement}
Data sharing not applicable to this article as no datasets were generated or analyzed during the current study.

\providecommand{\bysame}{\leavevmode\hbox to3em{\hrulefill}\thinspace}
\providecommand{\MR}{\relax\ifhmode\unskip\space\fi MR }
\providecommand{\MRhref}[2]{%
  \href{http://www.ams.org/mathscinet-getitem?mr=#1}{#2}
}
\providecommand{\href}[2]{#2}


\begin{thebibliography}{10}

\bibitem{Alves}
C.~J.~S. Alves, \emph{On the choice of source points in the method of fundamental solutions}, Eng. Anal. Bound. Elem. \textbf{33} (2009), 1348-1361.

\bibitem{Alves-Antunes_2005}
C.~J.~S. Alves and P.~R.~S. Antunes, \emph{The method of fundamental solutions
  applied to the calculation of eigenfrequencies and eigenmodes of 2{D} simply
  connected shapes}, Comput. Mater. Continua \textbf{2} (2005), 251--266.

\bibitem{Alves-Antunes_2013}
\bysame, \emph{The method of fundamental solutions applied to some inverse
  eigenproblems}, SIAM J. Sci. Comp. \textbf{35} (2013), A1689--A1708.

\bibitem{Alves-Leitao}
C.~J.~S. Alves and V.~M.~A. Leit\~{a}o, \emph{Crack analysis using an enriched MFS domain decomposition technique.}, Eng. Anal. Boundary Elem. \textbf{30}(3), (2006), 160--6.

\bibitem{Alves-Martins-Valtchev}
C.~J.~S. Alves, N.~F.~M. Martins and S.~S. Valtchev, \emph{Trefftz methods with cracklets and their relation to BEM and MFS.}, Eng. Anal. Boundary Elem. \textbf{95}, (2018), 93--104.

\bibitem{Alves-Silvestre}
C.~J.~S. Alves and A.~L. Silvestre, \emph{Density results using Stokeslets and a method of fundamental solutions for the Stokes equations.}, Eng. Anal. Boundary Elem., \textbf{28}(10)
(2004), 1245--1252.

\bibitem{Alves-Valtchev}
C.~J.~S. Alves and S.~S. Valtchev, \emph{A Kansa type method using fundamental solutions applied to elliptic PDEs.}, Advances in Meshfree Techniques (2007), 241--256.

\bibitem{Antunes_drum}
P.~R.~S. Antunes, \emph{Is it possible to tune a drum?}, J. Comput. Phys. \textbf{338}, (2017), 91--106.

\bibitem{MFSQR}
P.~R.~S. Antunes, \emph{Reducing the ill conditioning in the Method of Fundamental Solutions}, Adv. Comput. Math. \textbf{44}(1), (2018), 351--365.

\bibitem{Ant2}
P.~R.~S. Antunes, \emph{A numerical algorithm to reduce the ill-conditioning in meshless methods for the Helmholtz equation}, Numer. Algorithms \textbf{79}(3), (2018), 879--897.


\bibitem{Antunes_Valtchev}
P.~R.~S. Antunes and S.~S. Valtchev, \emph{A meshfree numerical method for acoustic wave propagation problems in planar domains with corner or cracks}, J. Comput. App. Math. \textbf{234}(9), (2010), 2646--2662.

\bibitem{Askour}
O. Askour, A. Tri, B. Braikat, H. Zahrouni and M. Potier-Ferry, \emph{Method of fundamental solutions and high order algorithm to solve nonlinear elastic problems.}, Eng. Anal. Boundary Elem., \textbf{89}
(2018), 25--35.


\bibitem{Barnett}
A.H. Barnett and T. Betcke, \emph{Stability and convergence of the method of
fundamental solutions for Helmholtz problems on analytic domains.}, J. Comput. Phys., \textbf{227}(14) (2008), 7003--26.

\bibitem{Berger}
J.~R. Berger and A. Karageorghis, \emph{The method of fundamental solutions for layered elastic materials.}, Eng. Anal. Boundary Elem. \textbf{25}(10), (2001); 877--86.

\bibitem{Bog}
A. Bogomolny, \emph{Fundamental solutions method for
elliptic boundary value problems}, SIAM J. Numer. Anal., 22(4)
(1985) pp.~644--669.

\bibitem{trefethen}
P.~D. Brubeck, Y. Nakatsukasa and L.~N. Trefethen, \emph{Vandermonde with Arnoldi}, to appear in SIAM Review.

\bibitem{CSChen}
C.~S. Chen, H. A. Cho and M. A. Golberg, \emph{Some comments on the ill-conditioning of the method of fundamental solutions}, Eng. Anal. Boundary Elem., \textbf{30}
(2006), 405--410.

\bibitem{chennh}
C.~S. Chen, C.~M. Fan and P.~H. Wen, \emph{The method of approximate particular solutions for solving certain partial differential equations.}, Numer. Methods Partial Differ. Eqs. \textbf{28}(2), (2012), 506--22.

\bibitem{Chen-Kara}
C.~S. Chen, A. Karageorghis ahd Y. Li, \emph{On choosing the location of the sources in the MFS.}, Numer. Algorithms \textbf{72}(1), (2016), 107--30.


\bibitem{Chen-electro}
C.~S. Chen, S.~Y. Reutskiy and V.~Y. Rozov, \emph{The method of fundamental solutions and its modifications for electromagnetic field problems.}, Comput. Assist. Mech. Eng. Sci., \textbf{16}(1) (2009), 21--33.

\bibitem{JTChen2}
J.T. Chen, J.-L. Yang, Y.-T. Lee and Y.-L. Chang, \emph{Formulation of the MFS for the two-dimensional Laplace equation with an added constant and constraint.}, Eng. Anal. Boundary Elem. \textbf{46}, (2014); 96--107.

\bibitem{JTChen}
J.T. Chen, C. S. Wu, Y. T. Lee and K. H. Chen, \emph{On the equivalence of the Trefftz method and method of fundamental solutions for Laplace and biharmonic equations.}, Comput. Math. Appl. \textbf{57} (2007), 851--79.

\bibitem{Cheng}
A.~H.~D. Cheng and Y. Hong, \emph{An overview of the method of fundamental solutions - Solvability, uniqueness, convergence and stability.}, Eng. Anal. Boundary Elem. \textbf{120}, (2020), 118--152.

\bibitem{zcli}
F. Dou,  L.-P. Zhang, Z.-C. Li and C.~S. Chen, \emph{Source nodes on elliptic pseudo-boundaries in the method of fundamental solutions for Laplace's equation; selection of pseudo-boundaries.}, J. Comput. Appl. Math. \textbf{337}, (2020); 112861.



\bibitem{FK}
G. Fairweather and A. Karageorghis, \emph{The method of fundamental solutions for elliptic boundary value problems}, Adv. Comput. Math. \textbf{9} (1998), 69--95.



\bibitem{Hon}
Y.C. Hon and M. Li, \emph{A discrepancy principle for the source points location in using the MFS for solving the BHCP.}, Int. J. Comput. Methods \textbf{6}, (2009), 181--197.

\bibitem{kar1} 
A. Karageorghis, \emph{A practical algorithm for determining the optimal pseudo-boundary in the method of fundamental solutions.}, Adv. Appl. Math. Mech. \textbf{1}(4), (2009), 510--28. 

\bibitem{kar-lesnic} 
A. Karageorghis and D. Lesnic, \emph{The method of fundamental solutions for the Oseen steady-state viscous flow past known or unknown shapes.}, Numer. Methods Partial Differ. Eqs. \textbf{35}(6), (2019), 2103--19.

\bibitem{Katsurada1}
M. Katsurada, \emph{A mathematical study of the charge simulation method. II}, J. Fac. Sci. Univ. Tokyo Sect. IA Math. \textbf{36} (1) (1989) 135--162.

\bibitem{Katsurada2}
\bysame, \emph{Charge simulation method using exterior mapping functions}, Jpn. J. Ind. Appl. Math. \textbf{11}(1),  (1994), 47--61.

\bibitem{Kitagawa}
T. Kitagawa, \emph{On the numerical stability of the method of fundamental
solution applied to the Dirichlet problem.}, Japan J. Appl. Math. \textbf{5} (1988), 123--33.

\bibitem{Kupradze}
V. D. Kupradze and M. A. Aleksidze, \emph{The method
of functional equations for the approximate solution of certain
boundary value problems}; U.S.S.R. Computational Mathematics and
Mathematical Physics \textbf{4} (1964), 82--126.

\bibitem{Li-CS-Kara}
M. Li, C.~S. Chen and A. Karageorghis, \emph{The MFS for the solution of harmonic boundary value problems with non-harmonic boundary conditions.}, Comput. Math. Appl. \textbf{66}(11), (2013);2400--24.

\bibitem{Sarler}
Q.~G. Liu and B. \u{S}arler, \emph{A non-singular method of fundamental solutions for the two-dimensional steady-state isotropic thermoelasticity problems.}, Eng. Anal. Boundary Elem. \textbf{75}, (2017);89--102.

\bibitem{mat-john}
R. Mathon and R. L. Johnston, \emph{The approximate solution of elliptic boundary-value problems by fundamental solutions.}, SIAM J. Numer. Anal. \textbf{14}, (1977), 638--50.

\bibitem{Martins-Silvestre}
N.~F.~M. Martins and A.~L. Silvestre, \emph{An iterative MFS approach for the detection of immersed obstacles.}, Eng. Anal. Boundary Elem., \textbf{32}(6)
(2008), 517--524.

\bibitem{raman} 
P. A. Ramachandran, \emph{Method of fundamental solutions: singular value decomposition analysis.}, Commun. Numer. Methods Eng. \textbf{18}, (2002), 789--801.

\bibitem{Rippa}
S. Rippa, \emph{An algorithm for selecting a good value for the parameter $c$ in radial basis function interpolation.}, Adv. Comput. Math. \textbf{11}, (1999), 193--210.


\bibitem{smy-kar} 
Y-S Smyrlis and A. Karageorghis, \emph{Some aspects of the method of fundamental solutions for certain harmonic problems.}, J. Sci. Comput. \textbf{16}(3), (2001), 341--71. 


\bibitem{Tsai}
C.~C. Tsai, D.~L. Young, C.~M. Fan and C.~W. Chen, \emph{MFS with time-dependent fundamental solutions for unsteady Stokes equations.}, Eng. Anal. Boundary Elem., \textbf{30}(10)
(2006), 897--908.

\bibitem{Young}
D.~L. Young and J.~W. Ruan, \emph{Method of fundamental solutions for scattering problems of electromagnetic waves.}, Comput. Model. Eng. Sci. \textbf{7}(2)
(2005), 223--32.




\end{thebibliography}
\end{document}